\documentclass[leqno]{article}
\usepackage{amssymb,amscd}
\usepackage{amsmath}
\usepackage{stmaryrd,mathabx}
\usepackage{rotating}
\usepackage[all]{xy}
\usepackage{caption,graphicx}

\def\un{{\rm 1\mkern-4mu I}}

\numberwithin{equation}{section}

\title{Geometry over ${\mathbb F}_1$}
\author{S. Haran \\ {\footnotesize (haran@tx.technion.ac.il)}}

\begin{document}

\maketitle

\vglue 2cm

\begin{enumerate}
\item[0.] Introduction

\item[1.] The ``Field with one element'' ${\mathbb F}$

\item[2.] $C{\mathbb F} R^t$: commutative ${\mathbb F}$-rings with involution

\item[3.] $C {\mathcal G} R$ : commutative generalized rings

\item[4.] Basic facts

\item[5.] Valuation generalized ring and the ``zeta machine''

\item[6.] Modules and differentials

\item[7.] ${\mathbb N}$ and ${\mathbb Z}$ as generalized rings

\item[8.] (Symmetric) Geometry

\item[9.] Pro-schemes

\item[10.] Arithmetical surface, and new (ordinary) commutative rings

\item[11.] Completed vector bundles
\end{enumerate}

\newpage

\setcounter{section}{-1}
\section{Introduction}

Andr\'e Weil, in a letter to his sister Simone Weil, compares the basic mystery of mathematics to the Roseta stone. But the Roseta stone has three languages (Demotic and Greek that were well known, and the Hieroglyphic that was a complete mystery), all describing one and the same reality, while in mathematics we have one language, addition and multiplication, the language of commutative rings, and this language is describing three similar, but different, realities: the mysterious number fields -- finite extensions of the rational numbers ${\mathbb Q}$; function fields -- finite extensions of ${\mathbb F}_q (t)$, which correspond to maps $f : X \to {\mathbb P}^1$ of smooth proper curves over ${\mathbb F}_q$; finite extensions of ${\mathbb C} (t)$, which correspond to maps $f : X \to {\mathbb P}^1$ of compact Riemann surfaces. 

\smallskip

When Alexander Grothendieck founded the modern language of algebraic geometry, he had behind him functional analysis, and the Gelfand theorem on the equivalence of the category of compact Hausdorff spaces and the opposite of the category of commutative $C^*$-algebras, and he had in front of him the Weil conjectures, which are all above ${\rm spec} \, {\mathbb Z}$, and so he choses the language of commutative rings as the foundations of geometry. Comparing arithmetic and geometry we see that this language is not general enough: the ``integers'' at a real (resp. complex) place of a number field ${\mathbb Z}_{\mathbb R}^1 = [-1,1]$ (resp. ${\mathbb Z}_{\mathbb C}^1 = \{ z \in {\mathbb C} , \vert z \vert \leq 1 \}$) are {\it not} closed under addition, and are not commutative rings; the ``arithmetical plane'' reduces to its diagonal ${\mathbb Z} \otimes {\mathbb Z} = {\mathbb Z}$, unlike the geometric analogue $F[t] \underset{F}{\otimes} F[t] = F[t_1 , t_2]$; we are missing the ``field with one element'' -- the common field of all the finite fields ${\mathbb Z} / (p)$, just as ${\mathbb F}_q$ is the common field of all the finite fields ${\mathbb F}_q [t] / (p(t))$, $p(t)$ irreducible.

\smallskip

We have some hints: the group ${\rm GL}_n ({\mathbb Q}_p)$ acts on ${\mathbb Q}_p^{\oplus n}$, and the stabilizer of the lattice
\begin{equation}
\label{eq1}
{\mathbb Z}^{\oplus n}_p = \left\{ (x_1 , \ldots , x_n) \in {\mathbb Q}_p^{\oplus n} \, , \quad \max \{ \vert x_1 \vert , \ldots , \vert x_n \vert \} \leq 1 \right\}
\end{equation}
is the maximal compact subgroup ${\rm GL}_n ({\mathbb Z}_p) \subseteq {\rm GL}_n ({\mathbb Q}_p)$. Similarly the group ${\rm GL}_n ({\mathbb R})$ (resp. ${\rm GL}_n ({\mathbb C})$) act on ${\mathbb R}^{\oplus n}$ (resp. ${\mathbb C}^{\oplus n}$), and the stabilizer of the unit-$L_2$-ball
\begin{equation}
\label{eq2}
\mbox{``} \, {\mathbb Z}_{\mathbb R}^n \, \mbox{''} = \left\{ (x_1 , \ldots , x_n) \in {\mathbb R}^{\oplus n} \, , \ \sum_{i \, = \, 1}^n \vert x_i \vert^2 \leq 1 \right\}
\end{equation}
resp.
$$
\mbox{``} \, {\mathbb Z}_{\mathbb C}^n \, \mbox{''} = \left\{ (x_1 , \ldots , x_n) \in {\mathbb C}^{\oplus n} \, , \ \sum_{i \, = \, 1}^n \vert x_i \vert^2 \leq 1 \right\}
$$
is the maximal compact subgroup $\mbox{``} \, {\rm GL}_n ({\mathbb Z}_{\mathbb R}) \, \mbox{''} = O(n) \subseteq {\rm GL}_n ({\mathbb R})$ the orthogonal group (resp. $\mbox{``} \, {\rm GL}_n ({\mathbb Z}_{\mathbb C}) \, \mbox{''} = U(n) \subseteq {\rm GL}_n ({\mathbb C})$ the unitary group). Another basic hint is the Cauchy-Schwartz inequality: for $x = (x_1 , \ldots , x_n)$, $y = (y_1 , \ldots , y_n)$ in ${\mathbb Z}_{\mathbb R}^n$ (resp. ${\mathbb Z}_{\mathbb C}^n$) $x \sslash y := x_1 y_1 + x_2 y_2 + \ldots + x_n y_n$ is in ${\mathbb Z}_{\mathbb R}^1$ (resp. ${\mathbb Z}_{\mathbb C}^1$). The essential thing in (\ref{eq2}) is that we use the $L_2$-norm (a fact that goes back to Pythagorass!), unlike the $L_{\infty}$-norm we have in (\ref{eq1}). In one dimension, there is no difference between $L_2$ and $L_{\infty}$ norms, and so we will have higher dimensions as part of our language.

\smallskip

We remark that the very same problem exists in physics: the interval of speeds $(-c,c)$, ($c$ being the speed of light in vacuum), is not closed under addition. Einstein's solution was to replace addition by the operation $\frac{x+y}{1 + \frac{xy}{c^2}}$. But this operation does not work for the complex numbers: the set $\{ x \in {\mathbb C} , \vert x \vert < c \}$ is not closed under this operation; it is closed under the operation $\frac{x+y}{1+ \frac{x \bar y}{c^2}}$, with $\bar y$ ($=$ conjugate of $y$), in the denominator, but this operation is not commutative or associative.

\smallskip

We shall replace Grothendieck's commutative rings, $C{\rm Ring}$, by the category of ``commutative generalized rings'' $C {\mathcal G} R$, and even by the more general category of ``commutative ${\mathbb F}$-rings with involution'', $C{\mathbb F} R^t$. These give a diagram of adjunctions:
\begin{equation}
\label{eq3}
\xymatrix{
&C {\mathcal G} R \ \ar@{^{(}->}@<1ex>[r]^{\mathcal F} \ar@<-1ex>[d]_{{\mathbb N} \, \underset{\mathbb F}{\otimes}} &C{\mathbb F} R^t \ar@{->>}@<1ex>[l]^{\mathcal U} \ar@<-1ex>[d]_{{\mathbb N} \, \underset{\mathbb F}{\otimes}} \\
C{\rm Rig} \ \ar@{^{(}->}@<1ex>[r] \ar@<-1ex>[d]_K &{\mathbb N} \backslash C {\mathcal G} R \ar@{_{(}->}@<-1ex>[u] \ar@{->>}@<1ex>[l] \ \ar@{^{(}->}@<1ex>[r]^{\mathcal F} \ar@<-1ex>[d]_{{\mathbb Z} \, \underset{\mathbb N}{\otimes}} &{\mathbb N} \backslash C{\mathbb F} R^t \ar@{_{(}->}@<-1ex>[u] \ar@{->>}@<1ex>[l]^{\mathcal U} \ar@<-1ex>[d]_{{\mathbb Z} \, \underset{\mathbb N}{\otimes}} \\
C{\rm Ring} \ \ar@{^{(}->}@<1ex>[r] \ar@{_{(}->}@<-1ex>[u] &{\mathbb Z} \backslash C {\mathcal G} R \ar@{_{(}->}@<-1ex>[u] \ar@{->>}@<1ex>[l] \ \ar@{^{(}->}@<1ex>[r]^{\mathcal F} &{\mathbb Z} \backslash C{\mathbb F} R^t \ar@{->>}@<1ex>[l]^{\mathcal U} \ar@{_{(}->}@<-1ex>[u]
}
\end{equation}
In this diagram the left adjoint ${\mathcal F}$ is written above the right adjoint ${\mathcal U}$, and we have ${\mathcal U} \circ {\mathcal F} \cong {\rm id}$. The category $C{\rm Rig}$ is the category of ``commutative rings without negatives'', i.e. $(R , + , \cdot , 0,1)$, $R$ is a set with two associative and commutative operations of addition ``$+$'' and multiplication ``$\cdot$'', with units $0$ and $1$, with multiplication distributive over addition, and with functions that preserve the operations and units.
\begin{equation}
\label{eq4}
\mbox{E.g.} \qquad R := \{0,1\}_{\max} \subseteq [0,1]_{\max} \subseteq [0,\infty)_{\max}
\end{equation}
with the operations
$$
\begin{matrix}
&x + y := \max \{ x,y \} \hfill \\
&x \cdot y := \mbox{ordinary multiplication.}
\end{matrix}
$$
For proofs of all the claims below we refer to \cite{13}.

\smallskip

This paper would have not existed without the encouragement of Laurent Lafforgue and the amazing typing job done by C\'ecile Gourgues and the artistic pictures done by Marie-Claude Vergne.

\section{The ``field with one element'' ${\mathbb F}$}

We capture the ``field with one element'' via the category of ``finite dimensional ${\mathbb F}_1$-vector spaces''. Denote by ${\mathbb F}$ the category with objects the finite sets, and with maps the partial-bijections:
\begin{equation}
\label{eq11}
{\mathbb F}_{Y,X} = {\mathbb F} (X,Y) = \left\{ X \supseteq D (\varphi) \underset{^\sim}{\overset{\varphi}{\longrightarrow}} I(\varphi) \subseteq Y \right\} .
\end{equation}

\bigskip

\noindent {\bf Remark:} Adding a ``zero'' to a finite set $X$, we obtain the pointed set $X_0 = X \amalg \{O_X \}$, and a map $\varphi \in {\mathbb F} (X,Y)$ gives a map $\varphi_0 \in {\rm Set}_0 (X_0 , Y_0)$ with $\varphi_0 (x) = \left\{ \begin{matrix} \varphi (x) &x \in D(\varphi) \\ O_Y \hfill &x \notin D(\varphi) \end{matrix} \right.$, we thus get an identification
\begin{equation}
\label{eq12}
{\mathbb F} (X,Y) \equiv \left\{ \varphi_0 \in {\rm Set}_0 (X_0 , Y_0), {\rm coker} \ker \varphi_0 \overset{_\sim}{\longrightarrow} \ker {\rm coker} \, \varphi_0 \right\} .
\end{equation}
The category ${\mathbb F}$ has kernels and cokernels, and in (\ref{eq12}) we see that they commute, which is one of the axioms of an abelian category.

\bigskip

\noindent {\bf Remark:} For $\varphi \in {\mathbb F} (X,Y)$ we get a $Y$ by $X$ matrix with values $0,1$:
\begin{equation}
\label{eq13}
(\widetilde\varphi_{y,x}) \in \{0,1\}^{Y \times X} \, , \quad \widetilde\varphi_{y,x} = \left\{ \begin{matrix} 1 &y = \varphi (x) \hfill \\ 0 &\mbox{otherwise.} \hfill \end{matrix} \right.
\end{equation}
We can thus identify ${\mathbb F} (X,Y)$ with the $\{0,1\}$-valued, $Y$ by $X$ matrices having at most one $1$ in every raw and in every column.

\smallskip

Note that ${\mathbb F}$ is self-dual: $( \ )^t : {\mathbb F} \overset{_\sim}{\longrightarrow} {\mathbb F}^{\rm op}$, where for
$$
\varphi = \left\{ D(\varphi) \underset{^\sim}{\overset{\varphi}{\longrightarrow}} I(\varphi) \right\} \in {\mathbb F} (X,Y) \, , \quad \varphi^t = \left\{ I(\varphi) \underset{^\sim}{\overset{\varphi^{-1}}{\longrightarrow}} D(\varphi) \right\} \in {\mathbb F} (Y,X) 
$$
(or taking the transpose of the $\{0,1\}$ matrix $\widetilde\varphi$):
\begin{equation}
\label{eq14}
(\varphi \circ \psi)^t = \psi^t \circ \varphi^t \, , \quad ({\rm id}_X)^t = {\rm id}_X \, , \quad \varphi^{tt} = \varphi \, .
\end{equation}
The category ${\mathbb F}$ does not have sums or products, but disjoint union induces a symmetric monoidal structure on ${\mathbb F}$:
\begin{equation}
\label{eq15}
\oplus : {\mathbb F} \times {\mathbb F} \to {\mathbb F} \, , \quad X \oplus Y = X \amalg Y \, .
\end{equation}
The unit is the empty set $[0]$, which is the initial and final object of ${\mathbb F}$. We have
\begin{equation}
\label{eq16}
(\varphi \oplus \psi)^t = \varphi^t \oplus \psi^t \, .
\end{equation}

\section{$C{\mathbb F} R^t$: commutative ${\mathbb F}$-rings with involution}

For $B \in C{\rm Rig}$ we associate the category ${\mathbb F} (B)$ of finitely generated free $B$-modules with a basis. Thus the objects are the finite sets, and for finite sets $X,Y$ the arrows from $X$ to $Y$ are given by the $B$-valued, $Y$ by $X$, matrices: $(B)^{Y \times X} \ni b = (b_{y,x})$. The operation of composition is matrix multiplication $\left[ (b'_{z,y}) \circ (b_{y,x}) \right]_{z,x} = \underset{y}{\sum} \ b'_{z,y} \cdot b_{y,x}$. This category is self-dual via transpose of matrices, and has a symmetric monoidal structure via direct-sum of matrices. Note that $(B)^{Y \times X}$ contains ${\mathbb F} (X,Y) \subseteq \{0,1\}^{Y \times X}$, and we obtain a functor ${\mathbb F} \hookrightarrow {\mathbb F}(B)$, which is the identity on objects, and is a strict-monoidal functor. This motivates the following

\bigskip

\noindent (2.1) \quad {\bf Definition.} An ${\mathbb F}$-ring $A$ is a category, with a functor $\varepsilon : {\mathbb F} \hookrightarrow A$, which is the identity on objects, and such that $\varepsilon [0]$ is the initial and final object of $A$, and we have a symmetric monoidal structure
$$
\oplus : A \times A \to A
$$
such that $\varepsilon$ is strictly monoidal: the associativity, commutativity, and unit isomorphisms are those of ${\mathbb F}$. 

\bigskip

We say $A$ is an ${\mathbb F}$-ring with involution, if it is self-dual via an involution that extends the involution on ${\mathbb F}$
$$
\xymatrix{
{\mathbb F} \ \ar@{^{(}->}[r] \ar[d]_{s} &A \ar[d]_{s}^{( \ )^t} \\
{\mathbb F}^{\rm op} \ \ar@{^{(}->}[r] &A^{\rm op}
}
$$
and is compatible with the symmetric monoidal structure
$$
(a_0 \oplus a_1)^t = a_0^t \oplus a_1^t \, .
$$
Thus the objects of $A$ are the finite sets, and we have to give the arrows from a finite set $X$ to a finite set $Y : A_{Y,X} = A(X,Y)$.

\smallskip

We have the following:
\begin{enumerate}
\item[(0)] $A_{Y , [0]} = \{0\} = A_{[0],X}$
\item[(1)] ${\mathbb F}_{Y,X} \subseteq A_{Y,X}$
\item[(2)] $\circ : A_{Z,Y} \times A_{Y,X} \to A_{Z,X}$ associative, unital, extends $\circ$ on ${\mathbb F}$
\item[(3)] $\oplus : A_{Y_0 , X_0} \times A_{Y_1 , X_1} \to A_{Y_0 \amalg Y_1 , X_0 \amalg X_1}$ with associativity, commutativity, and unit isomorphisms coming from ${\mathbb F}$
\item[(4)] $( \ )^t : A_{Y,X} \to A_{X,Y}$, $a^{tt} = a$, $(a_1 \circ a_2)^t = a_2^t \circ a_1^t$, $(a_0 \oplus a_1)^t = a_0^t \oplus a_1^t$.
\end{enumerate}

All the operations in (2), (3), (4) agree with the corresponding operations in ${\mathbb F}$. 

\smallskip

A homomorphism $\varphi \in {\mathbb F} R^t (A,B)$ is a functor over ${\mathbb F}$, strict-monoidal, and preserving the involution. Thus we have functions $\varphi = \varphi_{Y,X} : A_{Y,X} \to B_{Y,X}$ for $X,Y \in {\mathbb F}$, and we have
\begin{eqnarray}
\varphi (a_1 \circ a_2) &= &\varphi (a_1) \circ \varphi (a_2) \nonumber \\
\varphi (a_0 \oplus a_1) &= &\varphi (a_0) \oplus \varphi (a_1) \nonumber \\
\varphi (a^t) &= &\varphi (a)^t \nonumber \\
\varphi (\varepsilon_A (\psi)) &= &\varepsilon_B (\psi) \quad \mbox{for} \ \psi \in {\mathbb F}_{Y,X} \, . \nonumber 
\end{eqnarray}
Thus we have a category ${\mathbb F} R^t$.

\smallskip

Note that for $X,Y,Z \in {\mathbb F}$, $a \in A_{Y,X}$, we have $\underset{Z}{\oplus} \, a \in A_{\underset{Z}{\oplus} \, Y , \, \underset{Z}{\oplus} \, X} = A_{Z \times Y , \, Z \times X}$. We use this in the following critical

\bigskip

\noindent (2.2) \quad {\bf Definition.} An ${\mathbb F}$-ring (with involution) $A \in {\mathbb F} R^t$ is called {\it commutative} if for $X,Y,Z \in {\mathbb F}$, $a \in A_{Y,X}$, $b \in A_{[1], Z}$, $d \in A_{Z,[1]}$, we have in $A_{Y,X}$:
$$
a \circ \left[ \underset{X}{\oplus} \, (b \circ d) \right] = \left[ \underset{Y}{\oplus} \, (b \circ d) \right] \circ a = \left( \underset{Y}{\oplus} \, b \right) \circ \left( \underset{Z}{\oplus} \, a \right) \circ \left( \underset{X}{\oplus} \, d \right) .
$$
We denote by $C{\mathbb F} R^t$ the full subcategory of ${\mathbb F} R^t$ with objects  the commutative ${\mathbb F}$-rings with involution.

\smallskip

Note that for $B \in C{\rm Rig}$, we have ${\mathbb F} (B) \in C{\mathbb F} R^t$ is commutative, and we get a functor
$$
{\mathbb F} : C{\rm Rig} \hookrightarrow C{\mathbb F} R^t \, .
$$
It is a full and faithfull embedding: For $\varphi \in C{\mathbb F} R^t ({\mathbb F} (A) , {\mathbb F}(B))$, and $a = (a_{y,x}) \in {\mathbb F} (A)_{Y,X} = A^{Y \times X}$, since $\varphi$ is a functor over ${\mathbb F} : (\varphi_{Y,X} (a))_{y,x} = \varphi_{[1],[1]} (a_{y,x})$, and $\varphi$ is determined by $\varphi_{[1],[1]} : A \to B$ which preserves multiplication, and the units $0,1$, but also preserves addition:
\setcounter{equation}{2}
\begin{eqnarray}
\varphi_{[1],[1]} (a_1 + a_2) &= &\varphi \left( (1,1) \circ (a_1 \oplus a_2) \circ (1,1)^t \right) \\
&= &(\varphi (1) , \varphi (1)) \circ \left( \varphi (a_1) \oplus \varphi (a_2) \right) \circ (\varphi (1) , \varphi (1))^t  \nonumber \\
&= &(1,1) \circ \left( \varphi (a_1) \oplus \varphi (a_2) \right) \circ (1,1)^t \nonumber \\
&= &\varphi_{[1],[1]} \, (a_1) + \varphi_{[1],[1]} \, (a_2) \, . \nonumber
\end{eqnarray}

Note that the ${\mathbb F}$-ring $A = {\mathbb F} (B)$, $B \in C{\rm Rig}$, satisfies the extra ``total-commutativity'' of the following

\bigskip

\noindent (2.4) \quad {\bf Definition.} An ${\mathbb F}$-ring $A \in {\mathbb F} R^t$ is called ``total-commutative'' if for all $X,Y,Z,W \in {\mathbb F}$, $a \in A_{Y,X}$, $b \in A_{W,Z}$, we have in $A_{Y \times W , \, X \times Z}$:
$$
\left( \underset{W}{\oplus} \, a \right) \circ \left( \underset{X}{\oplus} \, b \right) = \left( \underset{Y}{\oplus} \, b \right) \circ \left( \underset{Z}{\oplus} \, a \right) .
$$
This common value can be denoted by $a \otimes b = (a \otimes {\rm id}_W) \circ ({\rm id}_X \otimes b) =  ({\rm id}_Y \otimes b) \circ (a \otimes {\rm id}_Z)$ and we get another symmetric monoidal structure $\otimes : A \times A \to A$, which is the usual product on objects $X \otimes Z = X \times Z$; the unit is the one point set $[1]$; and $\otimes$ is distributive over $\oplus$, so $A$ is a $C{\rm Rig}$-category, cf.~[H07].

\smallskip

As we shall see, total commutativity implies ${\mathbb Z} \otimes {\mathbb Z} = {\mathbb Z}$, and it is not needed for the developement of geometry, and we shall not assume it in general.

\bigskip

\noindent (2.5) \quad {\bf Remark.} For a (small) category $G$, and for $A \in {\mathbb F} R^t$, we have the category of functors $A^G$, and we let $R_A^+ (G)$ denote its isomorphism (i.e. natural equivalences) classes, and let $[f]$ denote the class of $f : G \to A$. We have the well defined operations on $R_A^+ (G)$

$$
[f_0] \oplus [f_1] \left( \underset{g_1}{\bullet} \overset{g}{\longleftarrow} \underset{g_0}{\bullet} \right) \equiv \left[ \underset{f_0 (g_1) \oplus f_1 (g_1)}{\bullet} \xleftarrow{ \ f_0 (g) \oplus f_1 (g) \ } \underset{f_0 (g_0) \oplus f_1 (g_1)}{\bullet} \right]
$$
$$
[f_0] \otimes [f_1] \left( \underset{g_1}{\bullet} \overset{g}{\longleftarrow} \underset{g_0}{\bullet} \right) \equiv \left[ \underset{f_0 (g_1) \times f_1 (g_1)}{\bullet} \xleftarrow{ \ \left( \underset{f_1 (g_1)}{\bigoplus} f_0 (g) \right) \circ \left(\underset{f_0 (g_0)}{\bigoplus} f_1 (g)\right) \ } \underset{f_0 (g_0) \times f_1 (g_0)}{\bullet} \right]
$$
$$
[f]^t \left( \underset{g_1}{\bullet} \overset{g}{\longleftarrow} \underset{g_0}{\bullet} \right) \equiv \left[ \underset{f(g_0)}{\bullet} \xleftarrow{ \ f(g)^t \ } \underset{f(g_1)}{\bullet} \right]
$$
making $R_A^+ (G)$ into a (non-commutative) rig with an involution. We let $R_A (G) = K R_A^+ (G)$ denote the associated ring. If $A$ is totally-commutative then $R_A^+ (G)$ and $R_A(G)$ are commutative. E.g. for a group $G$, and for $B \in C{\rm Ring},A=\mathbb{F}(B)$, $R_A(G)$ is the representation ring of $G$ on finitely-generated free $B$-modules, and $R_{\mathbb F} (G)$ is the Burnside ring of finite $G$-sets.

\bigskip

\noindent (2.6) \quad {\bf Remark.} If one wants to pass from ``finitely-generated-free-modules'' to ``finitely-generated-projective-modules'' (e.g. in the previous remark (2.5)), a quick way is to replace $A \in {\mathbb F} R^t$ by its idempotent completion $A^!$, with objects pairs $(X,p)$, $X \in {\mathbb F}$, $p \in A_{X,X}$, $p \circ p = p$, and mappings
$$
A^! ((X,p) , (Y,q)) = \{ f \in A_{Y,X} \, , \ f = q \circ f \circ p \} = q \circ A_{Y,X} \circ p \, .
$$

\section{$C {\mathcal G} R$: commutative generalized rings}

For $A \in C{\mathbb F} R^t$, we can forget its ``matrices'' $A_{Y,X} = A(X,Y)$, and just remember its ``vectors'' $A_X = A_{[1],X}$ ($\cong A_{X,[1]}$ since we have an involution), together with the operations of multiplication (resp. ``contraction'') by the direct sum of such vectors (resp. transposed). Thus for $B \in C{\rm Rig}$, let ${\mathcal G} (B) : {\mathbb F} \to {\rm Set}_0$, denote the functor
\begin{equation}
\label{eq31}
X \mapsto {\mathcal G} (B)_X := B^X \, .
\end{equation}
For a function of fintie sets $f : X \to Y$, put
$$
{\mathcal G} (B)_f = \prod_{y \, \in \, Y} {\mathcal G} (B)_{f^{-1} (y)} \quad (\equiv B^X) \, .
$$
Then we have the operations of

\bigskip

\noindent (3.2) \quad multiplication: ${\mathcal G}  (B)_Y \times {\mathcal G} (B)_f \to {\mathcal G} (B)_X$
$$
\begin{matrix}
b^Y , b^f &\mapsto &b^Y \lhd b^f \hfill \\
&&(b^Y \lhd b^f)_x = b_{f(x)}^Y \cdot b_x^f
\end{matrix}
$$

\noindent (3.3) \quad contraction: ${\mathcal G}  (B)_X \times {\mathcal G} (B)_f \to {\mathcal G} (B)_Y$
$$
\begin{matrix}
b^X , b^f &\mapsto &b^X \sslash b^f \hfill \\
&&(b^X \sslash b^f)_y = \underset{x \, \in \, f^{-1} (y)}{\sum} \, b_x^X \cdot b_x^f \, .
\end{matrix}
$$

This motivates the following

\bigskip

\noindent (3.4) \quad {\bf Definition.} A generalized ring $A$ is a functor ${\mathbb F} \to {\rm Set}_0$, $X \mapsto A_X$, and putting for $X,Y \in {\mathbb F}$, $f \in {\rm Set} (X,Y)$, $A_f = \underset{y \, \in \, Y}{\prod} \, A_{f^{-1} (y)}$, we have the operations

\bigskip

\noindent (3.5) \quad multiplication: $ \lhd : A_Y \times A_f \to A_X$

\bigskip

\noindent (3.6) \quad contraction: $\sslash : A_X \times A_f \to A_Y$.

\bigskip

These operations can be extended, fiber by fiber, to give for $g \in {\rm Set} (Y,Z)$

\bigskip

\noindent (3.7) \quad {\bf multiplication}: $\lhd : A_g \times A_f \to A_{g \circ f} \, , \quad (a^g \lhd a^f)_z = a_z^g \lhd (a_y^f)_{y \, \in g^{-1} (z)}$

\bigskip

\noindent (3.8) \quad {\bf contraction}: $\sslash : A_{g \circ f} \times A_f \to A_g \, , \quad (a^{g \circ f} \sslash a^f)_z = a_z^{g \circ f} \sslash (a_y^f)_{y \, \in g^{-1} (z)}$.

\bigskip

We require the following axioms:

\bigskip

\noindent (3.9) \quad {\bf Associativity}:  for $W \xleftarrow{ \ h \ } Z \xleftarrow{ \ g \ } Y \xleftarrow{ \ f \ }X$, $a_h \in A_h$, $a_g \in A_g$, $a_f \in A_f$, we have in $A_{h \circ g \circ f}$:
$$
a_h \lhd (a_g \lhd a_f) = (a_h \lhd a_g) \lhd a_f \, .
$$

\bigskip

\noindent (3.10) \quad {\bf Unit}:  We have $1 \in A_{[1]}$, and for any $a \in A_X$,
$$
1 \lhd a = a = a \lhd (1)_{x \, \in X} \, .
$$

I.e. $A$ with only the operation $\lhd$ is an ${\mathbb F}$-operad: an operad, with compatible $S_n$-action on $A_{[n]}$, together with $S_{[n]} \subseteq S_{[m]}$-covariant embeddings and projections $\xymatrix{A_{[n]} \ \ar@{^{(}->}@<0,2ex>[r]^{j} &A_{[m]} \ar@{->>}@<1ex>[l]^{ \ \pi}}$, $\pi \circ j = {\rm id}_{A_{[n]}}$, for $n < m$. 

\smallskip

The operation of contraction is dual to that of multiplication:

\bigskip

\noindent {\bf Duality}: For $W \xleftarrow{ \ h \ } Z \xleftarrow{ \ g \ } Y \xleftarrow{ \ f \ }X$, $a \in A_g$, $c \in A_f$, $d \in A_{h \circ g \circ f}$, we have in $A_h$:
\setcounter{equation}{10}
\begin{equation}
\label{eq311}
d \sslash (a \lhd c) = (d \sslash c) \sslash a
\end{equation}
and similarly, for $b \in A_{g \circ f}$, $c \in A_f$, $e \in A_{h \circ g}$, we have in $A_h$:
\begin{equation}
\label{eq312}
(e \lhd c) \sslash b = e \sslash (b \sslash c)
\end{equation}
$$
\xymatrix{
X \ar[rr]^c \ar[d]_-d \ar[rrd]^>>>b &&Y \ar[d]^-a \ar[lld]_>>>{e} \\
W &&Z \ar[ll]
}
$$

A homomorphism of generalized rings $\varphi \in {\mathcal G} R (A,B)$, is a natural transformation of functors, preserving the operations of multiplication and contraction, and the unit $1 : \{ \varphi_X : A_X \to B_X , X \in {\mathbb F}\}$, and
$$
\varphi (a \lhd a') = \varphi (a) \lhd \varphi (a') \, , \quad \varphi (a \sslash a') = \varphi (a) \sslash \varphi (a') \, , \quad \varphi (1) = 1 \, .
$$

Thus we have the category of generalized rings ${\mathcal G} R$.

\smallskip

We say that $A \in {\mathcal G} R$ is {\it commutative} if we have:
\begin{equation}
\label{eq313}
a \lhd (b \sslash b') = (a \lhd b) \sslash b'
\end{equation}
and if moreover, for $a \in A_X$, $b \in A_f$, $f : X \to Y$, $c \in A_g$, $g : Z \to Y$ we have in $A_Z$:
\begin{equation}
\label{eq314}
(a \sslash b) \lhd c = (a \lhd f^* c) \sslash g^* b 
\end{equation}
\begin{equation}
\label{eq315}
\xymatrix{
&X \ \underset{Y}{\prod} \ Z \ar[ld]_-{f^* c} \ar[rd]^-{g^* b} \\
a-X \ar[rd]_f^b &&Z \ar[ld]_c^g & 
\overset{\mbox{$(f^* c)_x := c_{f(x)}$}}{\mbox{$(g^* b)_z := b_{g(z)}$}} \\
&Y
} 
\end{equation}

We let $C{\mathcal G} R$ denote the full subcategory of ${\mathcal G} R$ with objects the commutative ${\mathcal G} R$. 

\smallskip

Note that for $B \in C{\rm Rig}$, ${\mathcal G} (B)$ is a commutative generalized ring, and we have a full and faithfull embedding
\begin{equation}
\label{eq316}
{\mathcal G} : C{\rm Rig} \hookrightarrow C{\mathcal G} R \, .
\end{equation}
Indeed, for $\varphi \in C{\mathcal G} R \, ({\mathcal G} (A) , {\mathcal G}(B))$, $\varphi_X (a)_x = \varphi_{[1]} (a_x)$, so $\varphi$ is determined by $\varphi_{[1]} : A \to B$, which is multiplicative, preserves $0$ and $1$, but also addition:
\begin{eqnarray}
\label{eq317}
\varphi_{[1]} (a_1 + a_2) &= &\varphi \left( ((1,1) \lhd (a_i)) \sslash (1,1) \right) \\
&= &((\varphi (1) , \varphi (1)) \lhd (\varphi (a_i))) \sslash (\varphi (1) , \varphi (1)) \nonumber \\
&= &((1,1) \lhd (\varphi (a_i))) \sslash (1,1) \nonumber \\
&= &\varphi_{[1]} (a_1) + \varphi_{[1]} (a_2) \, . \nonumber
\end{eqnarray}
We say that $A \in {\mathcal G} R$ is {\it ``totally-commutative''} if for $b \in A_f$, $c \in A_g$, in the notations of (\ref{eq315}), we have in $A_{X \, \underset{Y}{\prod} \, Z \, \to Y}$:
\begin{equation}
\label{eq318}
b \lhd f^* c = c \lhd g^* b \, .
\end{equation}
Note that for $B \in C{\rm Rig}$, the associated generalized ring ${\mathcal G} (B)$ is totally-commutative; but total-commutativity imply again ${\mathbb Z} \otimes {\mathbb Z} = {\mathbb Z}$, and we do not need to assume it for developing geometry.

\smallskip

Given $A\in C\mathbb{F}R^t$ we have $UA\in C\mathcal{G}R$, with
\begin{equation}
\begin{gathered}
(\mathcal{U} A)_X:=A_{[1],X}\\
a_Y \lhd (a_y^f):=a_Y\circ \left( \bigoplus_{y\in Y} a_y^f \right) \\
a_X \sslash (a_y^f):= a_X\circ \left( \bigoplus_{y\in Y} a_y^f \right)^t
\end{gathered}
\end{equation}
where $a_X\in A_{[1],X},a_Y\in A_{[1],Y},a_y^f\in A_{[1],f^{-1}(y)},f:X\to Y$. The left adjoint functor $\mathcal F:C\mathcal{G}R\to C\mathbb{F}R^t$ can be described as follows. For $B\in C\mathcal{G}R,\mathcal{F}(B)_{Y,X}=D(B)_{Y,X}/\approx$, where $D(B)_{Y,X}$ is the collection of data $\left( b_y^f,b_x^g \right)$, $b_y^f\in B_{f^{-1}(y)},b_x^g\in B_{g^{-1}(x)},f:Z\to Y,g:Z\to X$, and we divide it by the equivalence relation $\approx$ given by $\left(b_y^f,b_x^g\right)\approx \left(\bar{b}_x^{\bar{f}},\bar{b}_y^{\bar{g}} \right)$ iff for any $A\in C\mathbb{F}R^t$, any $\varphi\in C\mathcal{G}R(B,\mathcal{U}A)$, we have in $A_{Y,X}$:
\begin{equation}
\left(\bigoplus_{y\in Y}\varphi(b_y^f)\right)\circ \left(\bigoplus_{x\in X}\varphi(b_x^g)\right)^t= \left(\bigoplus_{y\in Y}\varphi(\bar{b}_y^{\bar{f}})\right)\circ \left(\bigoplus_{x\in X}\varphi(\bar{b}_x^{\bar{g}})\right)^t.
\end{equation}
The unit of adjunction $\varepsilon_X:B_X\to (\mathcal{UF}B)_X=\mathcal{F}(B)_{[1],X}$ is given by
\begin{equation}
\varepsilon_X(b)=(b, \{1_x \}_{x\in X}/\approx,\;f=c_X:X\to[1],g=\mathrm{id}_X:X\to X.
\end{equation}
The co-unit of adjunction $\varepsilon ^\ast_{Y,X}:\mathcal{F}(\mathcal{U}A)_{Y,X}\to A_{Y,X}$ is given by
\begin{equation}
\varepsilon^\ast_{Y,X}\left ( (b_y^f,b_x^g)/\approx \right)=\left(\bigoplus_{y\in Y}b_y^f\right)\circ \left(\bigoplus_{x\in X}b_x^g\right)^t.
\end{equation}
\section{Basic facts}

\noindent (4.1) \quad The categories $C{\mathcal G} R$ and $C{\mathbb F} R^t$ are complete and co-complete. 

\bigskip

Limits, and filtered co-limits are formed in Set ``pointwise'': for a functor $i \leadsto A_i$, from $I$ to $C{\mathcal G} R$ or $C{\mathbb F} R^t$,
$$
\left( \varprojlim_{i \, \in I} \, A_i \right)_{\!\!X} = \varprojlim_{i \, \in \, I} \, (A_i)_X \, ;
$$
$$
\left( \varinjlim_{i \, \in I} \, A_i \right)_{\!\!X} = \varinjlim_{i \, \in \, I} \, (A_i)_X \, , \quad \mbox{$I$ filtered: for $i_1 , i_2 \in I$ have $i_1 \to i$, $i_2 \to i$.}
$$
Moreover, $C{\mathcal G} R$ and $C{\mathbb F} R^t$ have push out diagrams
$$
\xymatrix{
A \ar@{.>}[r] &A \underset{c}{\otimes} B \\
C \ar[u] \ar[r] &B \ar@{.>}[u]
}
$$
and in the case of $C{\mathcal G} R$ some ``miracles'' happen that make $A \underset{c}{\otimes} B$ easier to describe (see below).

\bigskip

\noindent (4.2) \quad For $A \in C{\mathcal G} R$ (resp. $A \in C{\mathbb F} R^t$), we have a commutative, associative, unital monoid $A_{[1]}$ (resp. $A_{[1],[1]}$), with an involution $a^t := 1 \sslash a$ (resp. $a^t$), and it acts centrally on all the sets $A_X$ (resp. $A_{Y,X}$): For $a \in A_{[1]}$ (resp. $a \in A_{[1],[1]}$), $a_Y \in A_Y$ (resp. $a_{Y,X} \in A_{Y,X}$), $a \lhd a_Y = a_Y \lhd (a)_{y \in Y}$ (resp. $\left( \underset{Y}{\oplus} \, a \right) \circ a_{Y,X} = a_{Y,X} \circ \left( \underset{X}{\oplus} \, a \right)$). Moreover, if we have a homomorphism ${\mathbb N} \to A$ (resp. ${\mathbb Z} \to A$), then this is the multiplicative monoid of a $C{\rm Rig}$ (resp. $C{\rm Ring}$) with addition
$$
\begin{matrix}
a_1 + a_2 &:= &((1,1) \lhd (a_i)) \sslash (1,1) \hfill &A \in {\mathbb N} \backslash C{\mathcal G} R \hfill \\
&:= &(1,1) \circ (a_1 \oplus a_2) \circ (1,1)^t &A \in {\mathbb N} \backslash C{\mathbb F} R^t \nonumber
\end{matrix}
$$
 
\medskip

\noindent (4.3) \quad The initial object of $C{\mathcal G} R$ (resp. $C{\mathbb F} R$), is the ``field with one element'' ${\mathbb F}$, with ${\mathbb F}_X = X \amalg \{ 0_X \}$ (resp. ${\mathbb F}_{Y,X} = {\mathbb F} (X,Y)$).

\bigskip

\noindent (4.4) \quad There is a full and faithfull embedding of the category of commutative associative unital monoids $C{\rm Mon}$ into $C{\mathcal G} R$ (resp. $C{\mathbb F} R^t$), $M \mapsto {\mathbb F} \{M\}$ with ${\mathbb F} \{M\}_X = M \times X \amalg \{0_X \}$ (resp. ${\mathbb F} \{M\}_{Y,X} =$ set of $Y$ by $X$ matrices with values in $M \amalg \{0\}$ such that every raw/column contains at most one element of $M$).

\bigskip

\noindent (4.5) \quad There are free objects ${\mathbb F} \, [\delta_X] \in C{\mathcal G} R$ (resp. ${\mathbb F} \, [\delta_{Y,X}] \in C{\mathbb F} R^t$) such that for $A \in C{\mathcal G} R$ (resp. $C{\mathbb F} R^t$) we have:
$$
C{\mathcal G} R \, ({\mathbb F} \, [\delta_X] , A) \equiv A_X
$$
$$
\varphi \longleftrightarrow \varphi_X (\delta_X)
$$

$$
\mbox{resp.} \qquad C{\mathbb F} R^t \, ({\mathbb F} \, [\delta_{Y,X}] , A) \equiv A_{Y,X}
$$
$$
\qquad \qquad \qquad \varphi \longleftrightarrow \varphi_{Y,X} (\delta_{Y,X})
$$

\medskip

\noindent (4.6) \quad There are sub-$C{\mathcal G} R : {\mathbb Z}_{\mathbb R} \subseteq {\mathcal G} ({\mathbb R})$ and ${\mathbb Z}_{\mathbb C} \subseteq {\mathcal G} ({\mathbb C})$ (resp. sub-$C{\mathbb F} R^t : {\mathbb Z}_{\mathbb R} \subseteq {\mathbb F} ({\mathbb R})$ and ${\mathbb Z}_{\mathbb C} \subseteq {\mathbb F} ({\mathbb C})$) given by:
\begin{eqnarray}
({\mathbb Z}_{\mathbb R})_X &= &\left\{ a = (a_x) \in {\mathbb R}^X , \sum_{x \, \in \, X} \vert a_x \vert^2 \leq 1 \right\} \nonumber \\
({\mathbb Z}_{\mathbb C})_X &= &\left\{ a = (a_x) \in {\mathbb C}^X , \sum_{x \, \in \, X} \vert a_x \vert^2 \leq 1 \right\} \nonumber \\
({\mathbb Z}_{\mathbb R})_{Y,X} &= &\left\{ a = (a_{y,x}) \in {\mathbb R}^{Y \times X} , \ a (({\mathbb Z}_{\mathbb R})_X) \subseteq ({\mathbb Z}_{\mathbb R})_Y \right\} \nonumber \\
({\mathbb Z}_{\mathbb C})_{Y,X} &= &\left\{ a = (a_{y,x}) \in {\mathbb C}^{Y \times X} , \ a (({\mathbb Z}_{\mathbb C})_X) \subseteq ({\mathbb Z}_{\mathbb C})_Y \right\} . \nonumber
\end{eqnarray}
Moreover, there are ``residue-fields'' and surjective homomorphisms $\pi : {\mathbb Z}_{\mathbb R} \twoheadrightarrow {\mathbb F}_{\mathbb R}$, $\pi : {\mathbb Z}_{\mathbb C} \twoheadrightarrow {\mathbb F}_{\mathbb C}$, where
$$
({\mathbb F}_{\mathbb R})_X = \left\{ a = (a_x) \in {\mathbb R}^X , \sum_{x \, \in \, X} \vert a_x \vert^2 = 1 \right\} \amalg \{0_X\}
$$
are the real spheres augmented with a zero,
$$
({\mathbb F}_{\mathbb R})_{Y,X} = \left\{ {\mathbb R}^X \supseteq D(\varphi) \underset{^\sim}{\overset{\varphi}{\longrightarrow}}  I(\varphi) \subseteq {\mathbb R}^Y , \ \varphi \ \mbox{partial isometry} \right\}
$$
and similarly for ${\mathbb C}$. 

\smallskip

Note that for $A \in C{\mathcal G} R$ (resp. $A \in C{\mathbb F} R^t$) we have ``coefficient-map'':
$$
A_X \to (A_{[1]})^X \qquad (\mbox{resp.} \ A_{Y,X} \to (A_{[1],[1]})^{Y \times X})
$$
$$
a \mapsto (a \sslash 1_x)_{x \, \in \, X} \qquad (\mbox{resp.} \ a \mapsto (j_y^t \circ a \circ j_x)_{y \, \in \, Y , \, x \, \in \, X}).
$$
In most of the examples these are injections, but for ${\mathbb F}_{\mathbb R}$ and ${\mathbb F}_{\mathbb C}$ they are not!

\section{Valuation generalized ring and the ``zeta machine''}

For $K \in C{\mathcal G} R$ (resp. $C{\mathbb F} R^t$) we let $K^* = {\rm GL}_{[1]} (K)$ denote the invertible elements in $K_{[1]}$ (resp. $K_{[1],[1]}$). We say $K$ is a ``field'' if $K^* \amalg \{0\} = K_{[1]}$ (resp. $K_{[1],[1]}$); i.e. if every non-zero scalar is invertible. Let $B \subseteq K$ be a sub-$C{\mathcal G} R$ (resp. $C{\mathbb F} R^t$). It follows that $\{0\}$ is a prime of $B$, and we have the localization $B_{(0)} \subseteq K$ (see below). We say $B$ is {\it full} in $K$ if $B_{(0)} = K$. This means

\bigskip

\noindent (5.1) \quad $B \subseteq K$ {\it full} : For every $a \in K_X$ (resp. $K_{Y,X}$), there exist $d \in B_{[1]} \backslash \{0\}$, with $d \lhd a \in B_X$ 

\qquad (resp. $\left( \underset{Y}{\oplus} \, d \right) \circ a \in B_{Y,X}$).

\bigskip

The subset  of $K_X$ (resp. $K_{Y,X}$) consisting of elements $a$ such that for any $b \in B_X$ (resp. $b \in B_{X,[1]}$, $b' \in B_{[1],Y}$), $a \sslash b \in B_{[1]}$ (resp. $b' \circ a \circ b \in B_{[1],[1]}$) contains the set $B_X$ (resp. $B_{Y,X}$); if these sets are equal we say $B$ is {\it tame} in $K$:

\bigskip

\noindent (5.2) \quad $B \subseteq K$ {\it tame} : $B_X = \{ a \in K_X , B_X \sslash a \subseteq B_{[1]} \}$, 

\smallskip

\qquad  (resp. $B_{Y,X} = \{ a \in K_{Y,X} , B_{[1],Y} \circ a \circ B_{X,[1]} \subseteq B_{[1],[1]} \})$.

\bigskip

Finally, we say $B$ is a {\it valuation} sub-$C{\mathcal G} R$ (resp. $C{\mathbb F} R^t$) if it is full, and tame, and for any $a \in K^* : a \in B_{[1]}$ or $a^{-1} \in B_{[1]}$ (resp. $B_{[1],[1]}$). We let ${\rm Val} \, (K)$ denote all the valuation sub-$C{\mathcal G} R$ (resp. $C{\mathbb F} R^t$) of $K$. 

\smallskip

E.g. ${\mathcal G} ([0,1]) \subseteq {\mathcal G} ([0,\infty))$ is a valuation sub-$C{\mathcal G} R$; ${\mathbb F} ([0,1]) \subseteq {\mathbb F} ([0,\infty))$ a valuation sub-$C{\mathbb F} R^t$, cf. (0.4). 

\smallskip

We have the following interpretation of Ostrowski's theorem.

\bigskip

\noindent (5.3) \quad {\bf Ostrowski}: For a number field $K$,
$$
{\rm Val} \, (K) \equiv \{ K ; {\mathcal O}_{K,p} , \ p \subseteq {\mathcal O}_K \ \mbox{a finite prime}; \ K \cap \sigma^{-1} ({\mathbb Z}_{\mathbb C}) , \, \sigma : K \hookrightarrow {\mathbb C} \mod \sigma \sim \bar\sigma \} \, .
$$
E.g. for $K = {\mathbb Q}$,
$$
{\rm Val} \, ({\mathbb Q}) \equiv \{{\mathbb Q} \, ; \ {\mathbb Z}_{(p)} , \ p \ {\rm prime} ; \ {\mathbb Q} \cap {\mathbb Z}_{\mathbb R} \} \, .
$$

\bigskip

\noindent (5.4) \quad {\bf Remark}: Here ``valuation'' can be taken either in the sense of $C{\mathcal G} R$, or in the sense of $C{\mathbb F} R^t$, the results are the same. But note that if for $\sigma = 1/p \in [0,1]$, where $p \in [1,\infty]$, we let $({\mathbb Z}_{\mathbb R}^{(\sigma)})_{Y, \, X}$ denote the set of all real valued, $Y$ by $X$ matrices, that take the unit $L_p$-ball in ${\mathbb R}^X$ into the unit $L_p$-ball in ${\mathbb R}^Y$, then ${\mathbb Z}_{\mathbb R}^{(\sigma)}$ satisfy all the axioms of being a valuation sub-$C{\mathbb F} R$, but ${\mathbb Z}_{\mathbb R}^{(\sigma)}$ fails to have an involution, except in the $L_2$ case $\sigma = \frac12$.

\smallskip

In the above case of a number field $K$, the ordered abelian group $K^* / B^*$ is class 1, and can be viewed as an ordered subgroup of ${\mathbb R}^+ = (0,\infty)$. We get a homomorphism $\vert \ \vert_{[1]} : K^* \to (0,\infty)$, which we extend to
\setcounter{equation}{4}
\begin{equation}
\label{eq55}
\vert \ \vert_{[1]} : K_{[1]} \to [0,\infty) \, , \ \vert x \vert = 0 \Leftrightarrow x = 0 \, , \ \vert x_1 \circ x_2 \vert = \vert x_1 \vert \cdot \vert x_2 \vert \, .
\end{equation}
We also have 
\begin{equation}
\label{eq56}
\vert \ \vert_X : K_X \to [0,\infty) \, , \qquad (\mbox{resp.} \ \vert \ \vert_{Y, \, X} : K_{Y , \, X} \to [0,\infty))
\end{equation}
such that
$$
\vert a^Y \lhd a^f \vert_X \leq \vert a^Y \vert_Y \cdot \underset{y \, \in \, Y}{\rm Max} \, \vert a_y^f \vert_{f^{-1} (y)} \, , \  \vert a^X \sslash a^f \vert_Y \leq \vert a^X \vert_X \cdot \underset{y \, \in \, Y}{\rm Max} \, \vert a_y^f \vert_{f^{-1} (y)}
$$
$$
\mbox{(resp. $\vert a_1 \circ a_2 \vert \leq \vert a_1 \vert \cdot \vert a_2 \vert$, $\vert a_0 \oplus a_1 \vert \leq {\rm Max} \left\{ \vert a_0 \vert , \vert a_1 \vert \right\}$, $\vert a^t \vert = \vert a \vert$) }
$$
so that $B_X = \{ a \in K_X , \vert a \vert_X \leq 1 \}$, (resp. $B_{Y, \, X} = \{ a \in K_{Y, \, X} , \vert a \vert_{Y , \, X} \leq 1 \}$) is a sub-$C{\mathcal G} R$ (resp. $C{\mathbb F} R^t$) of $K$. 

\smallskip

Moreover, we have:
\begin{eqnarray}
\label{eq57}
\vert a \vert_X &= &{\rm inf} \left\{ \vert d \vert_{[1]}^{-1} \, , \ d \in K^* \, , \ d \lhd a \in B_X \right\} \quad \mbox{``full''} \\
&= &\sup \left\{ \vert b \sslash a \vert_{[1]} \, , \ b \in B_X \right\} \quad \mbox{``tame''} \nonumber
\end{eqnarray}
resp.
\begin{eqnarray}
\vert a \vert_{Y, \, X} &= &{\rm inf} \left\{ \vert d \vert_{[1]}^{-1} \, , \ d \in K^* \, , \ \left(\underset{Y}{\oplus} \ d \right) \circ a \in B_{Y, \, X} \right\} \quad \mbox{``full''} \nonumber \\
&= &\sup \left\{ \vert b' \circ a \circ b \vert_{[1]} \, , \ b' \in B_{[1],Y} \, , \ b \in B_{X,[1]} \right\} \quad \mbox{``tame''}.\nonumber
\end{eqnarray}

\bigskip

Let $B$ be a compact topological valuation $C{\mathbb F} R^t$. This means the sets $B_{Y, \, X}$ have a compact topology, such that the operations of composition, direct-sum, and transpose, are continuous. The field $K = B_{(0)}$ is then locally compact, and assume the valuation is given by a norm map
\begin{eqnarray}
\label{eq58}
\vert \ \vert_{Y, \, X} : K_{Y, \, X} &\to &[0,\infty) \\
\mbox{\begin{rotate}{90}$\subseteq$\end{rotate}} &&\quad\mbox{\begin{rotate}{90}$\subseteq$\end{rotate}} \nonumber \\
B_{Y, \, X} &\to &[0,1] \, . \nonumber
\end{eqnarray}
E.g. $B = {\mathbb F} ({\mathbb Z}_p)$ or $B = {\mathbb Z}_{\mathbb R}$.

\smallskip

We let
\begin{equation}
\label{eq59}
S_B^{[n]} = \{ b \in B_{[1],[n]} \, , \ \vert b \vert_{[1],[n]} = 1 \}
\end{equation}
denote the ``sphere'' in $n$-dimensional space, on which we have the homogeneous action of the compact group
\begin{equation}
\label{eq510}
{\rm GL}_{[n]} (B) = \{ b \in B_{[n],[n]} \, , \ \exists \, b' \in B_{[n],[n]} \, , \ b \circ b' = b' \circ b = {\rm id}_{[n]} \} \, .
\end{equation}
There exists a unique ${\rm GL}_{[n]} (B)$-invariant probability measure on $S_B^{[n]}$. We denote this Haar, or Maak, measure by $\sigma_B^{[n]}$. Assume we have a homomorphism ${\mathbb F} ({\mathbb N}) \to K$, so we have the elements
\begin{equation}
\label{eq511}
\un_{[n]} = (1,1,\ldots , 1) \in K_{[1],[n]} \, , \quad \un_{[n]}^t \in K_{[n],[1]} \, .
\end{equation}
Contracting $\sigma_B^{[n]}$ with $\un_{[n]}$ we get a probability measure on $K_{[1],[1]} = \{ 0 \} \amalg K^*$, and integrating $\vert \ \vert_{[1]}^{s-1}$ we get $[n]$-dimensional zeta
\begin{equation}
\label{eq512}
L_{[n]} (B,s) = \int_{S_B^{[n]}} \left\vert a \circ \un_{[n]}^t \right\vert_{[1]}^{s-1} \qquad \sigma_B^{[n]} (da) \, .
\end{equation}
Finally, we pass to the $n \to \infty$ limit to get the zeta of $B$:
\begin{equation}
\label{eq513}
L(B,s) = \lim_{n \to \infty} L_{[n]} (B,s) \, .
\end{equation}
E.g. We have the following calculations:
\begin{equation}
\label{eq514}
L ({\mathbb Z}_p , s) = \lim_{n \to \infty} \int_{{\rm Max} \, \{\vert a_1 \vert , \ldots , \vert a_n \vert\}=1} \vert a_1 + \ldots + a_n \vert_p^{s-1} \quad \sigma_{{\mathbb Z}_p}^{[n]} (da) = \frac{\zeta_p (s)}{\zeta_p (1)} \, , \ \zeta_p (s) = (1-p^{-s})^{-1} \, ,
\end{equation}
\begin{equation}
\label{eq515}
L ({\mathbb Z}_{\mathbb R} , s) = \lim_{n \to \infty} \int_{\vert a_1 \vert^2 + \ldots + \vert a_n \vert^2=1} \vert a_1 + \ldots + a_n \vert^{s-1} \quad \sigma_{{\mathbb Z}_{\mathbb R}}^{[n]} (da) = \frac{\zeta_{\mathbb R} (s)}{\zeta_{\mathbb R} (1)} \, , \ \zeta_{\mathbb R} (s) = 2^{\frac s2} \cdot \Gamma \left( \frac s2 \right) \, . 
\end{equation}

\section{Modules and differentials}

\noindent (6.1) \quad {\bf Definition.} For $A \in C{\mathbb F} R^t$, an $A$-module is a functor $M : A \times A^{\rm op} \to Ab$, $Y,X \mapsto M_{Y , \, X}$, such that $M_{[0], \, X} = \{0\} = M_{Y, \, [0]}$. Thus $M_{Y, \, X}$ are abelian groups, $Y,X \in {\mathbb F}$, and we have maps:
\begin{eqnarray}
A_{Y' , \, Y} \times M_{Y , \, X} \times A_{X, \, X'} &\to &M_{Y' , \, X'} \nonumber \\
a' \, , \ m \, , \ a &\mapsto &a' \circ m \circ a \nonumber
\end{eqnarray}
such that
\begin{eqnarray}
a' \circ (m_1 + m_2) \circ a &= &a' \circ m_1 \circ a + a' \circ m_2 \circ a \nonumber \\
a'_2 \circ (a'_1 \circ m \circ a_1) \circ a_2 &= &(a'_2 \circ a'_1) \circ m \circ (a_1 \circ a_2) \nonumber \\
{\rm id}_Y \circ m \circ {\rm id}_X&= &m \, . \nonumber
\end{eqnarray}
Such $A$-module is called ``commutative'' if for $m \in M_{Y, \, X}$, $b \in A_{[1] , \, Z}$, $d \in {\mathbb A}_{Z, \, [1]}$, we have:
$$
\left[ \underset{Y}{\oplus} \, (b \circ d) \right] \circ m \circ {\rm id}_X = {\rm id}_Y \circ m \circ \left[ \underset{X}{\oplus} \, (b \circ d) \right] = \left( \underset{Y}{\oplus} \, b \right) \circ \left( \underset{Z}{\oplus} \, m \right) \circ \left( \underset{X}{\oplus} \, d \right) .
$$

\bigskip

\noindent (6.2) \quad {\bf Definition.} For $A \in C{\mathcal G} R$, an $A$-module is a functor $M : {\mathbb F} \to Ab$, $X \mapsto M_X$, together with operations for $f : X \to Y$, $A_f = \underset{y \, \in \, Y}{\prod} \, A_{f^{-1} (y)}$,

\medskip

multiplication: $M_Y \times A_f \to M_X$, $m \lhd (a_y)$

\medskip

contraction: $M_X \times A_f \to M_Y$, $m \sslash (a_y)$,

\medskip

\noindent such that: $m \sslash f = m \lhd f^t = f(m)$ for $f \in {\mathbb F}_{Y, \, X}$, $m \in M_X$
\begin{eqnarray}
(m_1 + m_2) \lhd a &= &m_1 \lhd a + m_2 \lhd a \nonumber \\
(m_1 + m_2) \sslash a &= &m_1 \sslash a + m_2 \sslash a \nonumber \\
(m \lhd a_1) \lhd a_2 &= &m \lhd (a_1 \lhd a_2) \nonumber \\
(m \sslash a_1) \sslash a_2 &= &m \sslash (a_2 \lhd a_1) \nonumber \\
m \sslash (a_1 \sslash a_0) &= &(m \lhd a_0) \sslash a_1 \nonumber \\
m \lhd (a_1 \sslash a_0) &= &(m \lhd a_1) \sslash a_0 \nonumber 
\end{eqnarray}

and commutativity: $(m \sslash a) \lhd c = (m \lhd f^* c) \sslash g^* a$, $a \in A_f$, $c \in A_g$.

\bigskip

I.e. for $C{\mathbb F} R^t$ we take bi-modules (and we can add an involution on them that exchanges the left and right structures), while for $C {\mathcal G} R$ we take right modules. In both cases of $A \in C{\mathcal G} R$ and $A \in C{\mathbb F} R^t$, we get a complete and co-complete abelian category of $A$-modules with enough projective and injectives. A homomorphism $\varphi : A \to B$ gives an adjunction
\setcounter{equation}{2}
\begin{equation}
\label{eq63}
\xymatrix{
M^B &\mbox{$B$-mod} \ar@/^/[d] &N \ar@{|->}[d] \\
M \ar@{|->}[u] &\mbox{$A$-mod}\ar@/^/[u] &N_A
}
\end{equation}
For $A \in C{\mathbb F} R^t$, and $M \in A$-mod, we can form an ${\mathbb F}$-ring over $A$ (non-commutative, and without involution):
\begin{equation}
\label{eq64}
(A \, \pi M)_{Y , \, X} := A_{Y, \, X} \times M_{Y , \, X}
\end{equation}

\medskip

\noindent {\bf composition:} $(a_1 , m_1) \circ (a_2 , m_2) = (a_1 \circ a_2 , {\rm id}_Z \circ m_1 \circ a_2 + a_1 \circ m_2 \circ {\rm id}_X)$

\medskip

\noindent {\bf direct sum:} $(a_0 , m_0) \oplus (a_1 , m_1) = (a_0 \oplus a_1 , j_0^{Y_a} \circ m_0 \circ (j_0^X)^t + j_1^Y \circ m_1 \circ (j_1^X)^t)$ with $j_i^X \in {\mathbb F}_{X_0 \, \oplus \, X_1 , \, X_i}$ the inclusions.

\medskip

Similarly, for $A \in C{\mathcal G} R$, and $M \in A$-mod, we can form the (non-commutative) generalized ring over $A$:
\begin{equation}
\label{eq65}
(A \, \pi M)_X := A_X \times M_X
\end{equation}

\medskip

\noindent {\bf multiplication:} $(a_1 , m_1) \lhd (a_2 , m_2) = \left(a_1 \lhd a_2 , \ m_1 \lhd a_2 + \underset{y}{\sum} \, (m_2)_y \lhd (a_1\vert _y )_{x \, \in \, f^{-1} (y)} \right)$

\medskip

\noindent {\bf contraction:} $(a_1 , m_1) \sslash (a_2 , m_2) = \left( a_1 \sslash a_2 , \ m_1 \sslash a_2 + \underset{y}{\sum} \, (m_2)_y \sslash (a_1 \vert_{f^{-1} (y)}) \right)$.

\medskip

Fixing a homomorphism $k \to A$ of $C{\mathcal G} R$ or $C{\mathbb F} R^t$, we remember that $A \, \pi M$ is under $k$ via $k \to A \, \pi M$, $\lambda \mapsto (\lambda , 0)$, and we get adjunctions
\begin{equation}
\label{eq66}
\xymatrix{
C \ar@{|->}[d] &k \backslash {\mathcal G} R / A \ar@/^/[d] &k \backslash {\mathbb F} R / A \ar@/^/[d] &A \, \pi M \\
\Omega (C/k)^A  &\mbox{$A$-mod}\ar@/^/[u] &\mbox{$A$-mod} \ar@/^/[u] &M \ar@{|->}[u]
}
\end{equation}
for $C=A$:

$$
k \backslash {\mathcal G} R / A \, (A , A \, \pi M) \equiv {\rm Der}_k (A,M) \equiv \mbox{$A$-mod} \, (\Omega (A/k) , M)
$$
$$
\mbox{\hglue -20mm} \left( a \mapsto (a,\varphi (da)) \right) \leftarrow\!\!\!-\!\!\!-\!\!\!-\!\!\!-\!\!\!-\!\!\!-\!\!\!-\!\!\!-\!\!\!-\!\!\!-\!\!\!-\!\!\!-\!\!\!-\!\!\!-\!\!\!-\!\!\!-\!\!\!-\!| \ \varphi
$$
where $d : A \to \Omega (A/k)$ is the universal derivation (i.e. satisfying Leibnitz rule) trivial on $k$. From this all the usual properties of the Khaler differentials follow, e.g. for a triangle
$$
\xymatrix{
&B \\
k \ar[ur] \ar[r] &A \ar[u]
}
$$
we get an exact sequence of $B$-modules
\begin{equation}
\label{eq67}
\Omega (A/k)^B \to \Omega (B/k) \to \Omega (B/A) \to 0 \, .
\end{equation}
Moreover, the adjunction (\ref{eq66}) extends to a Quillen adjunction on the Quillen model categories of simplicial objects of the categories in (\ref{eq66}), and (\ref{eq67}) extends to a long exact sequence. The Quillen cotangent complex ${\mathbb L} \Omega (B/A)$ is an element of the derived category of $(\mbox{$B$-mod})^{\Delta^{\rm op}} \cong {\rm Ch}_{\bullet} (\mbox{$B$-mod})$.

\section{${\mathbb N}$ and ${\mathbb Z}$ as generalized rings}

We can now view the initial objects ${\mathbb N} \in C{\rm Rig}$, and ${\mathbb Z} \in C{\rm Ring}$, as $C{\mathcal G} R$, or as $C{\mathbb F} R^t$. Once we have the basic vector
$$
\delta = (1,1) \in {\mathbb N}^{[2]} = {\mathcal G} ({\mathbb N})_{[2]} = {\mathbb F} ({\mathbb N})_{[1],[2]}
$$
we have addition:
$$
a_1 + a_2 = (\delta \lhd (a_i)) \sslash \delta = \delta \circ (a_1 \oplus a_2) \circ \delta^t \, .
$$
It follows that $\delta$ generates ${\mathcal G} ({\mathbb N})$ and ${\mathbb F} ({\mathbb N})$. One can prove the relations are precisely:

\bigskip

\noindent (7.1) \quad $\delta$-{\bf unit}: $\delta \vert_i = 1$

\bigskip

\noindent (7.2) \quad $\delta$-{\bf commutativity}: $\delta \circ \begin{pmatrix} 0 &1 \\ 1 &0 \end{pmatrix} = \delta$

\bigskip

\noindent (7.3) \quad $\delta$-{\bf associativity}: $\delta \lhd (\delta ; 1) = \delta \lhd (1; \delta)$ or $\delta \circ (\delta \oplus 1) = \delta \circ (1 \oplus \delta)$
$$
{\includegraphics[width=50mm]{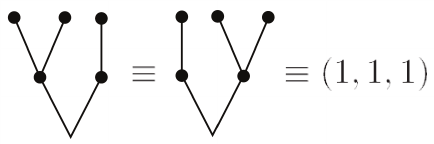}} 
$$
Moreover, ${\mathcal G} ({\mathbb Z})$ and ${\mathbb F} ({\mathbb Z})$, are generated by $\delta$ and by $(-1) \in {\mathcal G} ({\mathbb Z})_{[1]} = {\mathbb F} ({\mathbb Z})_{[1],[1]}$, with the above $\delta$-relations, and

\bigskip

\noindent (7.4) \quad $(-1) \circ (-1) = 1$

\bigskip

\noindent (7.5) \quad $\delta$-{\bf cancellation}: $(\delta \lhd (1;-1)) \sslash \delta = 0 = \delta \circ (1 \oplus (-1)) \circ \delta^t$.

\bigskip

Note that beside the above relations, we are assuming the relation of commutativity, and in particular:

\bigskip

\noindent (7.6) \quad $\delta^t \circ \delta = \begin{pmatrix} 1 \\ 1 \end{pmatrix} \circ (1,1) = \begin{pmatrix} 1 &1 \\ 1 &1 \end{pmatrix} = \begin{pmatrix} 1 &1 &0 &0 \\ 0 &0 &1 &1 \end{pmatrix} \circ \begin{pmatrix} 1 &0 \\ 0 &1 \\ 1 &0 \\ 0 &1 \end{pmatrix} = (\delta \oplus \delta) \circ (\delta \oplus \delta)^t$

\bigskip

or picturiously
$$
{\includegraphics[width=45mm]{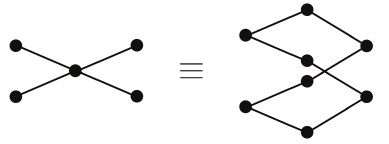}} 
$$

It follows that $d \, \delta$ generates $\Omega ({\mathbb N} / {\mathbb F})$ and $\Omega ({\mathbb Z} / {\mathbb F} \{ \pm 1 \})$, as ${\mathcal G} ({\mathbb N})$-module, and ${\mathcal G} ({\mathbb Z})$-module, respectively. In particular, for the $C {\mathcal G} R$ picture we get a map in degree $[1]$:

\bigskip

\noindent (7.7) \quad $d : {\mathbb N} \to \Omega ({\mathbb N} / {\mathbb F})_{[1]}$.

\bigskip

Here $\Omega ({\mathbb N} / {\mathbb F})_{[1]}$ turns out to be the free abelian group on generators
$$
\{ a,b \} := d \, \delta \sslash (a,b) \, , \qquad (a,b) \in {\mathbb N}^{[2]} \, ,
$$
modulo the relations:

\bigskip

\noindent (7.8) \quad $d \, (\delta\mbox{-unit}) : \{ a,0 \} = 0$

\bigskip

\noindent (7.9) \quad $d \, (\delta\mbox{-commutativity}) : \{ a,b \} = \{ b,a \}$

\bigskip

\noindent (7.10) \quad $d \, (\delta\mbox{-associativity}) : \{ a,b + c \} + \{ b,c \} = \{a+b,c\} + \{ a,b \}$

\bigskip

\noindent (7.11) \quad $d$ (commutativity).

\bigskip

The last relation (7.11) follows from, and is equivalent modulo 2-torsion to, the relation:

\bigskip

\noindent (7.12) \quad $\{ c \cdot a , c \cdot b \} = c \cdot \{ a , b \}$.

\bigskip

Denoting by $\Omega$ the quotient of $\Omega ({\mathbb N} / {\mathbb F})_{[1]}$ by the relation (7.12) and letting $\partial : {\mathbb N} \to \Omega$ denote $\frac12 \, d$ followed by the projection, we have:

\bigskip

\noindent (7.13) \quad Leibnitz: $\partial (n \cdot m) = n \, \partial (m) + m \, \partial (n)$

\bigskip

\noindent (7.14) \quad Almost-additivity: $\partial (n + m) = \partial (n) + \partial (m) + \{ n,m \}$

\bigskip

\noindent (7.15) \quad $\partial (0) = \partial (1) = 0$, and for $n > 1 : \partial (n) = \{1,1\} + \{1,2\} + \ldots + \{1,n-1\}$.

\bigskip

Note that it follows from Leibnitz that we have
\setcounter{equation}{15}
\begin{equation}
\label{eq716}
\partial (n) = \sum_p v_p (n) \cdot \frac np \, \partial (p)
\end{equation}
and $\Omega$ is the free abelian group on $\partial (p)$, $p$ prime ($v_p (n)$ the $p$-adic valuation).

\bigskip

\noindent (7.17) \quad {\bf Remark.} We can now also give the proof that total-commutativity imply ${\mathbb Z} \otimes {\mathbb Z} = {\mathbb Z}$ (a proof due to N.~Durov). We write on the left (resp. right) the formulas in the language of $C {\mathcal G} R$ (resp. $C {\mathbb F} R^t$). We note that ${\mathbb Z} \underset{{\mathbb F} \{ \pm 1 \}}{\otimes} {\mathbb Z}$ is generated (over ${\mathbb F} \{\pm 1\}$) by two generators $\delta_1$ and $\delta_2$ of degree $[2]$. Total-commutativity gives:
$$
\delta_1 \lhd (\delta_2 , \delta_2) = \delta_2 \lhd (\delta_1 , \delta_1) \qquad\qquad\qquad\qquad \delta_1 \circ (\delta_2 \oplus \delta_2) = \delta_2 \circ (\delta_1 \oplus \delta_1) 
$$

$$
{\includegraphics[width=120mm]{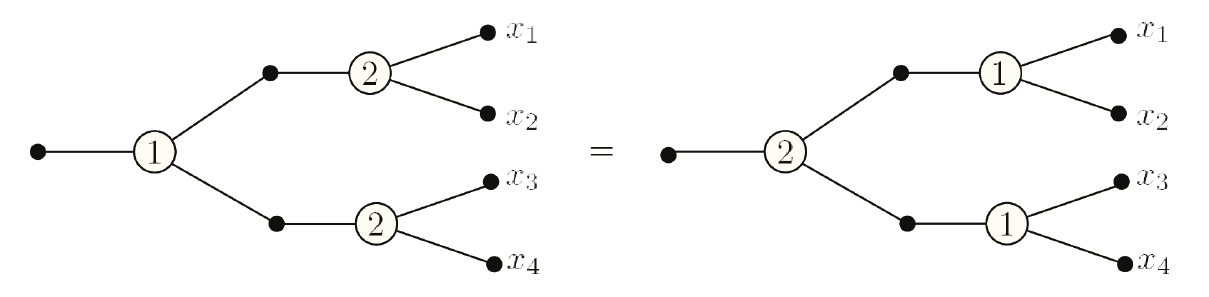}} 
$$
Restricting this identity to $\{ x_2 , x_3 \}$ gives:
$$
\delta_1 \lhd \left( \delta_2 \vert_{x_2} , \delta_2 \vert_{x_3} \right) = \delta_2 \lhd \left( \delta_1 \vert_{x_2} , \delta_1 \vert_{x_3} \right) \qquad \qquad \qquad \delta_1 \circ \left( \delta_2 \vert_{x_2} \oplus \delta_2 \vert_{x_3} \right) = \delta_2 \circ \left( \delta_1 \vert_{x_2} \oplus \delta_1 \vert_{x_3} \right)
$$

$$
{\includegraphics[width=110mm]{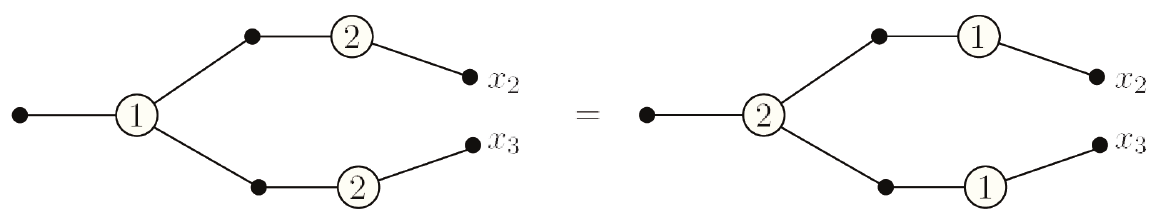}} 
$$
Using the $\delta$-unit relation, $\delta \vert_{x_i} = 1$, we get:
$$
{\includegraphics[width=60mm]{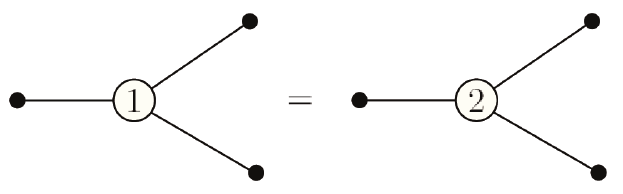}} 
$$
i.e. $\delta_1 = \delta_2$, and ${\mathbb Z} \underset{{\mathbb F} \{\pm 1\}}{\otimes} {\mathbb Z} = {\mathbb Z}$.

\section{(Symmetric) Geometry}

For $A \in C{\mathbb F} R^t$, or $A \in C {\mathcal G} R$, an {\it equivalence ideal} is a sub-$C{\mathbb F} R^t$, or $C {\mathcal G} R$, of the product ${\mathcal E} \subseteq A \, \Pi \, A$, which is an equivalence relation on $A$. We can then form the quotient $A/{\mathcal E}$, and we have a canonical homomorphism $\pi : A \twoheadrightarrow A/{\mathcal E}$. Given any homomorphism $\varphi : A \to B$, it gives rise to an equivalence ideal of $A$
$$
{\rm KER} (\varphi) = A \, \underset{B}{\Pi} \, A = \{ (a,a') , \varphi (a) = \varphi (a') \}
$$
and we have the factorization of $\varphi$
\begin{equation}
\label{eq81}
\xymatrix{
A \ar[r]^{\varphi} \ar@{->>}[d]_-{\pi} &B \\
A/{\rm KER} (\varphi) \ar[r]^-{\sim} &\varphi (A) \ar@{_{(}->}[u]
}
\end{equation}
We let ${\mathcal E}\!q (A)$ denote the equivalence ideals of $A$.

\bigskip

\noindent (8.2) \quad {\bf Definition:} For $A \in C {\mathcal G} R$ (resp. $C{\mathbb F} R^t$), an {\it ideal} of $A$ is a subset ${\mathfrak a} \subseteq A_{[1]}$ (resp. $A_{[1],[1]}$), such that for $X \in {\mathbb F}$, $a = (a_x) \in (a)^X$, $b,d \in A_X$ (resp. $b \in A_{[1], X}$, $d \in A_{X , [1]}$),
$$
(b \lhd (a_x)) \sslash d \in {\mathfrak a} \, , \quad \mbox{(resp.} \ b \circ \left( \underset{x \, \in \, X}{\oplus} \ a_x \right) \circ d \in {\mathfrak a}).
$$
We let $I(A)$ denote the set of ideals of $A$. There is a Galois correspondence:
\setcounter{equation}{2}
\begin{equation}
\label{eq83}
\xymatrix{
{\mathcal E}\!q (A) \ar@/^/[r] &I(A) \ar@/^/[l] 
}
\end{equation}
\begin{eqnarray}
{\mathcal E}_{\mathfrak a} \ \mbox{gen. by} \ (a,0) , \ a \in {\mathfrak a} &\leftarrow\!\!\!-\!\!\!-\!\!\!-\!\!\!-\!\!\!-\!\!\!-\!\!\!-\!\!\!-\!\!\!-\!\!\!-\!\!\!-\!| &{\mathfrak a} \nonumber \\
{\mathcal E} &|\!\!\!-\!\!\!-\!\!\!-\!\!\!-\!\!\!-\!\!\!-\!\!\!-\!\!\!-\!\!\!-\!\!\!-\!\!\!-\!\!\!\rightarrow &{\mathfrak a}_{\mathcal E} = \{ a , (a,0) \in {\mathcal E}_{[1]} \} \, . \nonumber
\end{eqnarray}

To motivate the following definition, imagine we would like to develope algebraic geometry based on commutative rings with a (possibly non-trivial) involution (cf. Remark 5.4). For such a ring $A$, we can forget the involution and look at ${\rm spec} \, (A)$: it is a locally ringed space with a (global) involution, but once we localize we loose the involution. A better way to proceed is to look at the subring of symmetric elements $A^+ = \{ a \in A , a^t = a \}$, and over the topological space ${\rm spec} \, (A^+)$ we have a sheaf of commutative rings with involution, with stalk over ${\mathfrak p} \in {\rm spec} \, (A^+)$ given by the localization $(A^+ \backslash \, {\mathfrak p})^{-1} \cdot A$.

\bigskip

\noindent (8.4) \quad {\bf Definition:} An ideal ${\mathfrak a}$ is called {\it symmetric} if it is generated as an ideal by its symmetric elements ${\mathfrak a}^+ = \{ a \in {\mathfrak a} , a^t = a \}$, i.e. for any $a \in {\mathfrak a}$, have $X \in {\mathbb F}$, $c = (c_x) \in ({\mathfrak a}^+)^X$, $c_x^t = c_x$, $b,d \in A_X$, with $a = (b \lhd c) \sslash d$ (resp. $b \circ \left( \underset{x \, \in \, X}{\oplus} \ c_x \right) \circ d^t$ for $C{\mathbb F} R^t$).

\bigskip

An ideal ${\mathfrak p} \subseteq A_{[1]}$ is called {\it prime} if $S_{\mathfrak p} = A_{[1]} \backslash \, {\mathfrak p}$ is {\it multiplicative}:
$$
1 \in S_{\mathfrak p} \, , \quad S_{\mathfrak p} \circ S_{\mathfrak p} = S_{\mathfrak p} \, .
$$
A symmetric-ideal ${\mathfrak p} \subseteq A_{[1]}$ is called {\it symmetric-prime} if $S_{\mathfrak p}^+ = A^+_{[1]} \backslash \, {\mathfrak p}$ is multiplicative.

\smallskip

We let $I^+ (A)$, ${\rm spec} \, (A)$, ${\rm spec}^+ (A)$ denote the set of symmetric ideals, primes, and symmetric-primes respectively.

\smallskip

The set ${\rm spec} \, (A)$ (resp. ${\rm spec}^+ (A)$) is not empty: every maximal proper (symmetric) ideal $m \in I(A)$ (resp. $m \in I^+ (A)$), which exists by Zorn, is a (symmetric) prime. The set ${\rm spec} \, (A)$ (resp. ${\rm spec}^+ (A)$) carry a topology, which is compact (every open cover has a finite sub-cover), and Zariski (every closed irreducible subset is the closure of a unique ``generic'' point), where the closed sets are
\setcounter{equation}{4}
\begin{equation}
\label{eq85}
V({\mathfrak a}) = \{ {\mathfrak p} \in {\rm spec} \, (A) , {\mathfrak p} \supseteq {\mathfrak a}\} \quad \mbox{(resp.} \ V^+ ({\mathfrak a}) = \{{\mathfrak p} \in {\rm spec}^+ (A) , {\mathfrak p} \supseteq {\mathfrak a}\})
\end{equation}
and a basis for the open sets is given by the ``basic''
\begin{equation}
\label{eq86}
D(f) = \{{\mathfrak p} \in {\rm spec} \, (A) , {\mathfrak p} \not\ni f \} , f \in A_{[1]} \, , \quad \mbox{(resp.} \ D^+ (f) = \{ {\mathfrak p} \in {\rm spec}^+ (A) , {\mathfrak p} \not\ni f \} , f = f^t \in A_{[1]}^+).
\end{equation}
On ${\rm spec} \, (A)$ we have a (global) involution ${\mathfrak p} \mapsto {\mathfrak p}^t$. There is a continuous map
\begin{eqnarray}
\label{eq87}
\pi : {\rm spec} \, (A) &\xrightarrow{ \ \ \ \ \ \ \ \ \ \ \ \ \ \ \ \ \ } &{\rm spec}^+ (A) \\
{\mathfrak p} &|\!\!\!-\!\!\!-\!\!\!-\!\!\!-\!\!\!-\!\!\!-\!\!\!-\!\!\!-\!\!\!-\!\!\!-\!\!\!-\!\!\!\rightarrow &\pi ({\mathfrak p}) = \mbox{ideal gen. by} \ {\mathfrak p}^+ \, . \nonumber
\end{eqnarray}
The process of localization with respect to a multiplicative set $S$ of symmetric elements $S \subseteq A_{[1]}^+$ works as usual for any $A \in C{\mathbb F} R^t$ or $A \in C {\mathcal G} R$, because the commutative monoid $A_{[1],[1]}^+$ or $A_{[1]}^+$ acts centrally on $A$. We obtain a homomorphism $\phi_S : A \to S^{-1} A$, which is universal taking $S$ to units. For an open subset ${\mathcal U} \subseteq {\rm spec}^+ (A)$, we have a $C {\mathcal G} R$ (resp. $C{\mathbb F} R^t$) ${\mathcal O}_A ({\mathcal U})$, and restriction homomorphisms $\rho_{\mathcal V}^{\mathcal U} : {\mathcal O}_A ({\mathcal U}) \to {\mathcal O}_A ({\mathcal V})$ for ${\mathcal V} \subseteq {\mathcal U}$, with $\rho_{\mathcal W}^{\mathcal V} \circ \rho_{\mathcal V}^{\mathcal U} = \rho_{\mathcal W}^{\mathcal U}$, $\rho_{\mathcal U}^{\mathcal U} = {\rm id}_{{\mathcal O}_A ({\mathcal U})}$, and for $X \in {\mathbb F}$ (resp. $Y,X \in {\mathbb F}$) ${\mathcal U} \mapsto {\mathcal O}_A ({\mathcal U})_X$ (resp. ${\mathcal O}_A ({\mathcal U})_{Y,X}$) is a sheaf of sets over ${\rm spec}^+ (A)$. The stalk at a point ${\mathfrak p} \in {\rm spec}^+ (A)$ is the localization
\begin{equation}
\label{eq88}
{\mathcal O}_A \vert_{\mathfrak p} = \varinjlim_{{\mathfrak p} \, \in \, {\mathcal U}} \, {\mathcal O}_A ({\mathcal U}) \xrightarrow{ \ \sim \ } (S_{\mathfrak p}^+)^{-1} A = A_{({\mathfrak p})}
\end{equation}
and it is ``local'' in the sense that it contains a unique maximal-symmetric-ideal $m_{\mathfrak p} \subseteq A_{({\mathfrak p})}$. The usual argument (using commutativity!) gives that for $f=f^t \in A_{[1]}^+$ we have over the basic open set:
\begin{equation}
\label{eq89}
{\mathcal O}_A (D^+ (f)) \xleftarrow{ \ \sim \ } (f^{\mathbb N})^{-1} A \, .
\end{equation}
All this is functorial: for a homomorphism $\varphi : A \to B$ we get a square of continuous maps
\begin{equation}
\label{eq810}
\xymatrix{
{\rm spec} \, (A) \ar[d]_{\pi_A} &&{\rm spec} \, (B) \ar[ll]_{\varphi_{[1]}^{-1}} \ar[d]^{\pi_B} \\
{\rm spec}^+ (A) &&{\rm spec}^+ (B) \ar[ll]_{\varphi^*} \\
}
\end{equation}
$$
\mbox{\hglue -48mm} \varphi^* (q) = \mbox{ideal gen. by} \ A_{[1]}^+ \cap \varphi^{-1} (q) \leftarrow\!\!\!-\!\!\!-\!\!\!-\!\!\!-\!\!\!-\!\!\!-\!\!\!-\!\!\!-\!\!\!-\!\!\!-\!\!\!-\!| \ q
$$
This leads to the introduction of the categories loc $C{\mathcal G} R$-Sp and loc $C{\mathbb F} R^t$-Sp, with objects $(X , {\mathcal O}_X)$, \qquad \quad $X \in {\rm Top}$ a topological space, ${\mathcal O}_X$ a sheaf of $C{\mathcal G} R$ (resp. $C{\mathbb F} R^t$) over $X$, with local stalks $m_x \subseteq ({\mathcal O}_X \vert_x)_{[1]}$ a unique maximal-symmetric-ideal, and mappings $(f,f^{\#}) : (X , {\mathcal O}_X) \to (Y , {\mathcal O}_Y)$ are pairs of $f \in {\rm Top} \, (X,Y)$ a continuous map, and $f^{\#} : {\mathcal O}_Y \to f_* \, {\mathcal O}_X$ a map of sheaves over $Y$ of $C{\mathcal G} R$ (resp. $C{\mathbb F} R^t$) and such that for $x \in X$ the composition
\begin{equation}
\label{eq811}
f_x^{\#} : {\mathcal O}_Y \vert_{f(x)} \xrightarrow{ \ f^{\#} \vert_{f(x)} \ } f_* \, {\mathcal O}_X \vert_{f(x)} = \varinjlim_{f(x) \, \in \, {\mathcal V}} {\mathcal O}_X (f^{-1} \, {\mathcal V}) \to \varinjlim_{x \, \in \, {\mathcal U}} {\mathcal O}_X ({\mathcal U}) = {\mathcal O}_X \vert_x
\end{equation}
is a {\it local-homomorphism}:
$$
f_x^{\#} (m_{f(x)}) \subseteq m_x \Leftrightarrow m_{f(x)} = (f_x^{\#})^* (m_x) \, .
$$

Within the category loc $C{\mathcal G} R$-Sp (resp. loc $C{\mathbb F} R^t$-Sp) we have the full subcategory $C{\mathcal G} R$-Sch (resp. $C{\mathbb F} R^t$-Sch) of ``schemes'' consisting of $(X , {\mathcal O}_X)$ which are locally affine: there is a covering by open sets \qquad \quad \qquad $X = \underset{i \, \in \, I}{\cup} \ {\mathcal U}_i$, and the canonical mappings give isomorphisms in loc $C{\mathcal G} R$-Sp (resp. loc $C{\mathbb F} R^t$-Sp) 

\noindent $({\mathcal U}_i , {\mathcal O}_X \vert_{{\mathcal U}_i}) \xrightarrow{ \ \sim \ } {\rm spec}^+ ({\mathcal O}_X ({\mathcal U}_i))$.

\smallskip

These categories $C{\mathcal G} R$-Sch and $C{\mathbb F} R^t$-Sch are much like ordinary schemes Sch: we can glue together objects $X_i$, along consistent glueing data on open subsets $X_i \supseteq {\mathcal U}_{ij} \underset{\sim}{\overset{\varphi_{ij}}{-\!\!\!-\!\!\!-\!\!\!\longrightarrow}} \, {\mathcal U}_{ji} \subseteq X_j$, $\varphi_{jk} \circ \varphi_{ij} = \varphi_{ik}$, $\varphi_{ii} = {\rm id}_{{\mathcal U}_i}$. Finite limits exists in $C{\mathcal G} R$-Sch and $C{\mathbb F} R^t$-Sch, in particular we have fiber products
\begin{equation}
\label{eq812}
\xymatrix{
&X \, \underset{Z}{\Pi} \, Y \ar[ld] \ar[rd] \\
X\ar[rd] &&Y \ar[ld] \\
&Z
}
\end{equation}

All the local calculations in $C{\mathcal G} R$-Sch (resp. $C{\mathbb F} R^t$-Sch) reduce to ``algebra'' in $C{\mathcal G} R$ (resp. $C{\mathbb F} R^t$), and the fiber product (\ref{eq812}) is obtained by glueing the affine pieces ${\rm spec}^+ \left( {\mathcal O}_X ({\mathcal U}) \underset{{\mathcal O}_Z({\mathcal W})}{\otimes} {\mathcal O}_Y ({\mathcal V}) \right)$.

\smallskip

Note that for $A \in C{\mathbb F} R^t$, the underlying topological space ${\rm spec}^+ (A)$ depends only on $UA \in C{\mathcal G} R$, and ${\rm spec}^+ (A) \equiv {\rm spec}^+ (UA)$, so that $C{\mathcal G} R$ is the basic geometric language. For $A \in C{\mathbb F} R^t$ we have more structure: the underlying symmetric monoidal category give rise to the algebraic $K$-theory spectrum:
\begin{equation}
\label{eq813}
K(A) \equiv {\mathbb Z} \times B \left( \underset{n \to \infty}{\lim} {\rm GL}_n (A)\right)^+ \equiv B (A^{-1} A) \equiv B\Gamma A(S^*)
\end{equation}
obtain by either Quillen $+$-construction, or ``group-localization'' $A^{-1} A$, or its delooping via Segal's $\Gamma$-category $\Gamma A$ applied to the sphere spectrum $S^*$.

\smallskip

Thus over any $C{\mathbb F} R^t$-Sch $X$ we have a canonical sheaf of symmetric spectra $K({\mathcal O}_X)$, and the associated sheaves of abelian groups $K_i ({\mathcal O}_X) = \pi_i \, K ({\mathcal O}_X)$.

\section{Pro-schemes}

The categories $C{\mathcal G} R$-Sch and $C{\mathbb F} R^t$-Sch are still not enough for arithmetic, and we pass to the categories of pro-objects: pro-$C{\mathcal G} R$-Sch and pro-$C{\mathbb F} R^t$-Sch. The objects are arbitrary functors $X : I \to C{\mathcal G} R$-Sch (resp. $C{\mathbb F} R^t$-Sch), $i \mapsto X_i$, $j \geq i \mapsto X_j \to X_i$, where we take for $I$ a countable partially ordered set, which is

\bigskip

\noindent (9.1) \quad {\bf directed:} for any $i_1 , i_2 \in I$, there is $j \in I$, $j \geq i_1$, $j \geq i_2$, and which is

\bigskip

\noindent (9.2) \quad {\bf cofinite:} for $j \in I$, $\{i \in I , j \geq i \}$ is finite.

\bigskip

The mappings from $\{X_i \}_{i \, \in \, I}$ to $\{Y_j \}_{j \, \in \, J}$ are given by
\setcounter{equation}{2}
\begin{equation}
\label{eq93}
\mbox{pro-$C{\mathcal G} R$-Sch} \, ( \{ X_i \}_{i \, \in \, I} , \{Y_j \}_{j \, \in \, J}) = \varprojlim_J \varinjlim_I \mbox{$C{\mathcal G} R$-Sch} \, (X_i , Y_j)
\end{equation}
and similarly with $C{\mathbb F} R^t$-Sch.

\smallskip

To motivate the passage to pro-objects, one may think of ordinary schemes, Sch, which is the full subcategory of loc Ring Sp, consisting of objects which are locally affine $\simeq {\rm spec} \, (A)$, $A \in C {\rm Ring}$. Given a functor $X : I \to {\rm loc \, Ring \, Sp}$, with $I$ directed as above, the limit $\underset{I}{\varprojlim} \, X_i$ exists in loc Ring Sp, and moreover, if $X_i = {\rm spec} \, (A_i)$ are all affine $\underset{I}{\varprojlim} \, {\rm spec} \, (A_i) = {\rm spec} \left( \underset{I}{\varinjlim} \, A_i \right)$ is again affine. But if $X_i \in {\rm Sch}$ are schemes for all $i \in I$, the limit $\underset{I}{\varprojlim} \, X_i$ need not be a scheme.

\smallskip

What happens is that for $X : I \to {\rm Sch} \subseteq \mbox{loc Ring Sp}$, a point $x = \{ x_i \}_{i \, \in \, I}$ in the space $\underset{I}{\varprojlim} \, X_i$ is a coherent family of points $x_i \in X_i$, and while each $x_i$ will have an affine neighborhood, these can shrink, so that $x$ will not have an affine neighborhood. The same phenomenon happens with $C{\mathcal G} R$-Sch and $C{\mathbb F} R^t$-Sch, which are not closed under directed inverse limits. The real and complex primes of a number field are such points that do not have an affine neighborhood.

\smallskip

Inside pro-$C{\mathcal G} R$-Sch, and pro-$C{\mathbb F} R^t$-Sch, we have the ``compactified'' $\overline{{\rm spec} \, {\mathbb Z}} = \{ X_N \}_{N \, \in \, I}$, with $I$ the set of square free integers, $N = p_1 \cdot p_2 \ldots p_{\ell}$ ($p_j$ distinct primes), with $M \geq N$ iff $N$ divide $M$, and with
\begin{equation}
\label{eq94}
X_N = {\rm spec} \, ({\mathbb Z}) \coprod_{{\rm spec} \, \left( {\mathbb Z} \left[ \frac1N \right]\right)} {\rm spec} \left( {\mathbb Z} \left[ \frac1N \right] \cap {\mathbb Z}_{\mathbb R} \right)
\end{equation}
i.e. we view ${\mathbb Z}$, ${\mathbb Z} \left[ \frac1N \right]$, $A_N = {\mathbb Z} \left[ \frac1N \right] \cap {\mathbb Z}_{\mathbb R}$ as $C{\mathcal G} R$ (resp. $C{\mathbb F} R^t$), the involution is the identity so the spectrum is the same as the symmetric spectrum, and we glue ${\rm spec} \, ({\mathbb Z})$ and ${\rm spec} \, (A_N)$ along the common basic open set
$$
{\rm spec} \, ({\mathbb Z}) \supseteq D_{\mathbb Z} (N) = {\rm spec} \left( {\mathbb Z} \left[ \frac1N \right]\right) = D_{A_N} \left( \frac1N \right) \subseteq {\rm spec} \, (A_N) \, .
$$
For $N$ dividing $M$, we have a natural map $X_M \xrightarrow{ \ f_N^M \ } X_N$, it is the identity on points, and moreover ${\mathcal O}_{X_N} \xrightarrow{ \ \sim \ } (f_N^M)_* \, {\mathcal O}_{X_M}$ is an isomorphism, but there are more open neighborhoods of the ``real-prime'' in $X_M$ than in $X_N$. Here the ``real-prime'' is the maximal ideal $\eta_N$ of $A_N = {\mathbb Z} \left[ \frac1N \right] \cap \, {\mathbb Z}_{\mathbb R}$, and for $N = p_1 \cdot p_2 \ldots p_{\ell}$, the specialization picture of the points of $X_N$ looks like:
\begin{equation}
\label{eq95}
\end{equation}
$$
{\includegraphics[width=130mm]{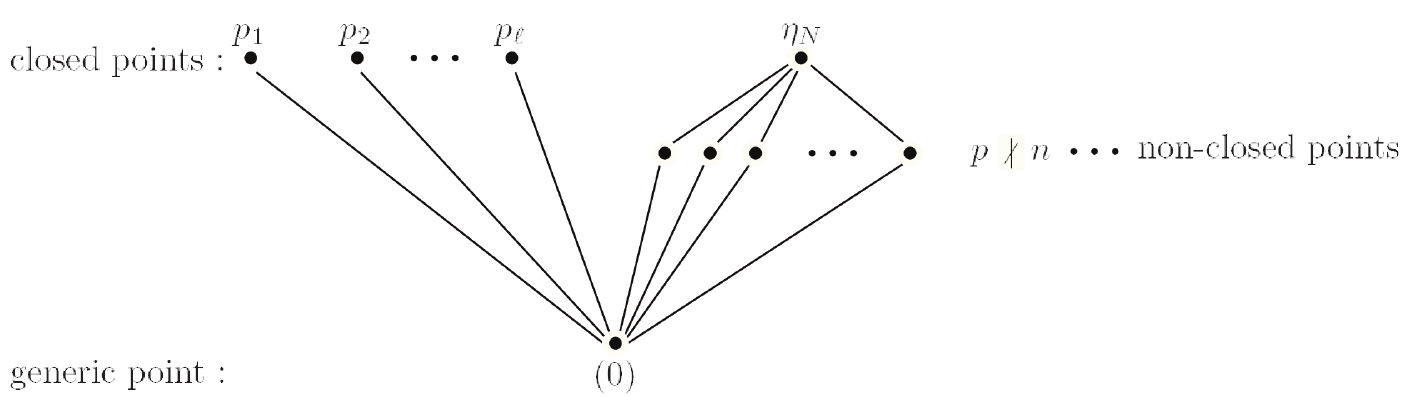}} 
$$

Note that $X_N$ looks like it has Krull dimension 2, but the inverse limit $X = \underset{I}{\varprojlim} \, X_N$ in loc -$C{\mathcal G} R$-Sp \quad (or in loc-$C{\mathbb F} R^t$-Sp) has Krull dimension 1, and the real-prime is a closed point just like the finite primes:

\bigskip

\noindent (9.6) \quad $X = \underset{I}{\varprojlim} \ X_N$: 
$$
{\includegraphics[width=60mm]{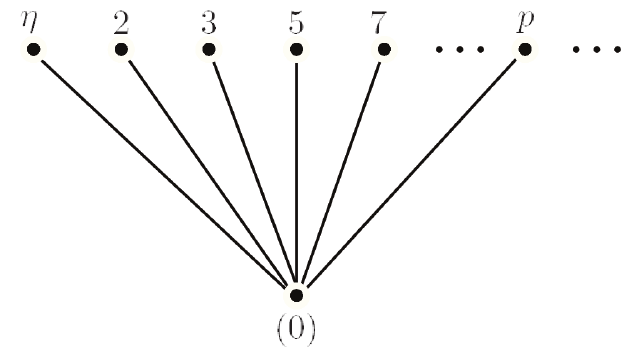}} 
$$

The open subsets of $X$ are arbitrary co-finite subsets ${\mathcal U} \subseteq X$, with $(0) \in {\mathcal U}$, and $X$ is integral:
\begin{equation}
\setcounter{equation}{7}
\label{eq97}
{\mathcal O}_X ({\mathcal U}) = {\mathbb Q} \cap \prod_{p \, \in \, {\mathcal U}} {\mathbb Z}_{(p)} \, , \quad \mbox{with} \quad {\mathbb Z}_{(\eta)} = {\mathbb Q} \cap {\mathbb Z}_{\mathbb R} \, .
\end{equation}
E.g. the global sections are
$$
{\mathcal O}_X (X) = {\mathbb Q} \cap \prod_{p \, \in \, X} {\mathbb Z}_{(p)} = {\mathbb Z} \cap {\mathbb Z}_{\mathbb R} = {\mathbb F} \, \{ \pm 1 \} \, .
$$

Similarly, for any number field $K$, with ring of integers ${\mathcal O}_K$, and with ``real and complex primes'' \qquad \quad ${\mathcal O}_{\sigma} = K \cap \sigma^{-1} ({\mathbb Z}_{\mathbb C})$, $\sigma \in \gamma_K = \{ \sigma : K \hookrightarrow {\mathbb C} \} / \sigma \sim \overline\sigma$, we have the compactified $\overline{{\rm spec} \, {\mathcal O}_K} = \{ X_N^K \}_{N \, \in \, I}$, with the same indexing set $I$ of square free integers $N = p_1 \ldots  p_{\ell}$ under divisibility, where $X_N$ is obtained by glueing ${\rm spec} \, ({\mathcal O}_K)$, and the various ${\rm spec} \left( {\mathcal O}_K \left[ \frac1N \right] \cap {\mathcal O}_{\sigma} \right)$, $\sigma \in \gamma_K$, along the common basic open set ${\rm spec} \left( {\mathcal O}_X \left[ \frac1N \right]\right)$. This gives $\overline{{\rm spec} \, {\mathcal O}_K}$ as an object of pro-$C{\mathcal G} R$-Sch or as an object of pro-$C{\mathbb F} R^t$-Sch. The global sections are ${\mathcal O}_{\overline{{\rm spec} \, {\mathcal O}_K}} \, (\overline{{\rm spec} \, {\mathcal O}_K}) = {\mathbb F} \{\mu_K\}$, where $\mu_K$ is the group of roots of unity in $K$. For $K \subseteq L$, there is a natural map $X_N^L \to X_N^K$ for all $N \in I$, and we get a map of pro-objects $\overline{{\rm spec} \, {\mathcal O}_L} \to \overline{{\rm spec} \, {\mathcal O}_K}$. Note that $\overline{{\rm spec} \, {\mathcal O}_K}$ is ``proper'' in the sense that the points of $\underset{N}{\varprojlim} \, X_N^K$ correspond bijectively with ${\rm Val} (K)$. Note that any $f \in K$ defines a geometric map from $\overline{{\rm spec} \, {\mathcal O}_K}$ to
\begin{equation}
\label{eq98}
{\mathbb P}^1 / {\mathbb F} = {\rm spec} \, ({\mathbb F} \{ z^{\mathbb N} \}) \coprod_{{\rm spec} \, ({\mathbb F} \{ z^{\mathbb Z}\})} {\rm spec} \left( {\mathbb F} \left\{ \left( \frac1z \right)^{\mathbb N} \right\}\right)
\end{equation}
via $z \mapsto f$, (and a ``Hurwitz-genus-formula'' for this map should give the $ABC$ conjecture \cite{23}).

\section{Arithmetical surface, and new commutative rings}

We have now new geometric objects, such as the Arithmetical-plane associated with a number field $K$
\begin{equation}
\label{eq101}
\xymatrix{
\overline{{\rm spec} \, {\mathcal O}_K} \underset{{\rm spec} \, {\mathbb F} \{\mu_K \}}{\prod} \overline{{\rm spec} \, {\mathcal O}_K} \ar[d] &= &\left\{ X_N^K \underset{{\rm spec} \, {\mathbb F} \{\mu_K \}}{\prod} X_M^K \right\}_{N, \, M \, \mbox{\footnotesize square-free}}  \\
\overline{{\rm spec} \, {\mathbb Z}} \underset{{\rm spec} \, {\mathbb F} \{\pm 1 \}}{\prod} \overline{{\rm spec} \, {\mathbb Z}} &= &\left\{ X_N \underset{{\rm spec} \, {\mathbb F} \{\pm 1 \}}{\prod} X_M \right\}_{(N,M) \, \in \, I \times I}
}
\end{equation}
and higher dimensional planes obtained by taking $n$-fold products. Note that the pro-objects in (\ref{eq101}) contain the ordinary affine open and dense sub-scheme
\begin{equation}
\label{eq102}
\xymatrix{
{\rm spec} \left( {\mathcal O}_K \underset{{\mathbb F} \{\mu_K \}}{\bigotimes} {\mathcal O}_K \right) \ar[d] \\
{\rm spec} \left( {\mathbb Z} \underset{{\mathbb F} \, \{ \pm 1 \}}{\bigotimes} {\mathbb Z} \right)
}
\end{equation}

When we have a homomorphism of $C{\mathcal G} R$ (resp. $C{\mathbb F} R^t$) $\varphi : {\mathbb Z} \to A$, we can make the commutative monoid $A_{[1]}$ (resp. $A_{[1],[1]}$) into an ordinary commutative ring, with an involution, by defining addition by
\begin{equation}
\label{eq103}
a_1 + a_2 := (\varphi ((1,1)) \lhd (a_i)) \sslash \varphi ((1,1))
\end{equation}
$$
\mbox{resp.} \ := \varphi ((1,1)) \circ (a_1 \oplus a_2) \circ \varphi ((1,1))^t
$$
(note that our axiom of commutativity imply the distributive law!).

\smallskip

We thus have the adjunction
\begin{equation}
\label{eq104}
C {\rm Ring} \, (B,A_{[1]}) \equiv {\mathbb Z} \backslash C{\mathcal G} R \ ({\mathcal G} (B),A)
\end{equation}
$$
C {\rm Ring} \, (B,A_{[1],[1]}) \equiv {\mathbb Z} \backslash C{\mathbb F} R^t \ ({\mathbb F} (B),A) \, .
$$
The unit of adjunction ${\mathcal G} (A_{[1]}) \to A$ (resp. ${\mathbb F} (A_{[1],[1]}) \to A$), give for any $A$ under ${\mathbb Z}$ an injection of ${\rm spec} \, (A)$ into the spectrum of an ordinary ring:
\begin{equation}
\label{eq105}
{\rm spec} \, (A) \hookrightarrow {\rm spec} \, (A_{[1]}) \, , \quad \mbox{(resp. ${\rm spec} \, (A_{[1],[1]})$)}.
\end{equation}
As a $C{\mathcal G} R$, or as $C{\mathbb F} R^t$, the integers ${\mathbb Z}$ are generated by the vector $\delta = (1,1)$, and the scalar $(-1)$. It follows that the $n$-fold tensor product
\begin{equation}
\label{eq106}
{\mathbb Z}^{\otimes n} = {\mathbb Z} \underset{{\mathbb F} \{\pm1\}}{\otimes} \ldots \underset{{\mathbb F} \{\pm1\}}{\otimes} {\mathbb Z}
\end{equation}
is generated by the vectors $\delta_1 , \ldots , \delta_n \in ({\mathbb Z}^{\otimes n})_{[2]}$, and $(-1)$. Thus every element of $({\mathbb Z}^{\otimes n})_{Y,X}$ can be represented as a graph made of the basic graphs
\begin{equation}
\label{eq107}
\end{equation}
$$
{\includegraphics[width=90mm]{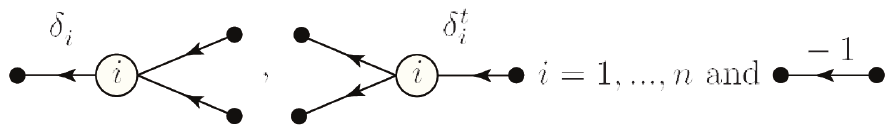}}
$$
and ``going from $X$ to $Y$''. Moreover, using the $\delta$-associativity relation (7.3), we can represent any $m$-fold product of a generator $\delta_i$ as the graph
\begin{equation}
\label{eq108}
\end{equation}
$$
{\includegraphics[width=65mm]{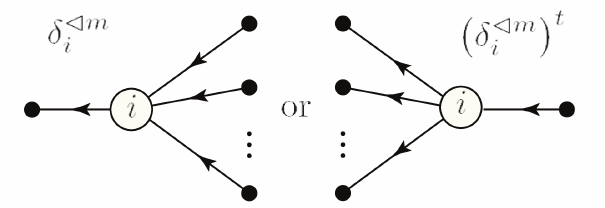}}
$$
Finally, we can use commutativity, to move all generators $\delta_i$ to the left, and all $\delta_j^t$ to the right. It follows that any element $a \in ({\mathbb Z}^{\otimes n})_{Y, \, X}$ has a representation as $a = (F_y , G_x , \sigma , \varepsilon)$, where $F_y$, $y \in Y$, $G_x$, $x \in X$, are finite-rooted-trees, together with a labeling of the non-leaves vertices
\begin{equation}
\label{eq109}
\coprod_{y \, \in \, Y} (F_y \backslash \partial F_y) \amalg \coprod_{x \, \in \, X} (G_x \backslash \partial G_x) \to \{ 1,2,\ldots , n\} \, ,
\end{equation}
and a bijection of this leaves
\begin{equation}
\label{eq1010}
\sigma :  \coprod_{y \, \in \, Y}\partial F_y \xrightarrow{ \ \sim \ } \coprod_{x \, \in \, X} \partial G_x \, ,
\end{equation}
and finally the signs
\begin{equation}
\label{eq1011}
\varepsilon :  \coprod_{y \, \in \, Y}\partial F_y \to \{\pm 1\} \, .
\end{equation}

\noindent (10.12) \qquad E.g.
$$
{\includegraphics[width=80mm]{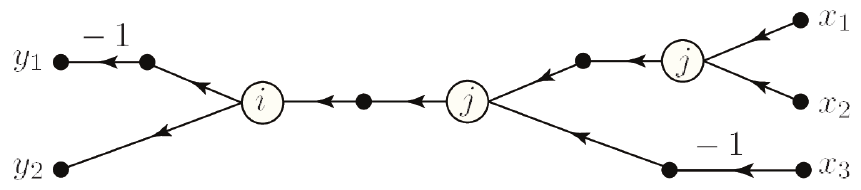}}
$$

\noindent moves to the equivalent
$$
{\includegraphics[width=90mm]{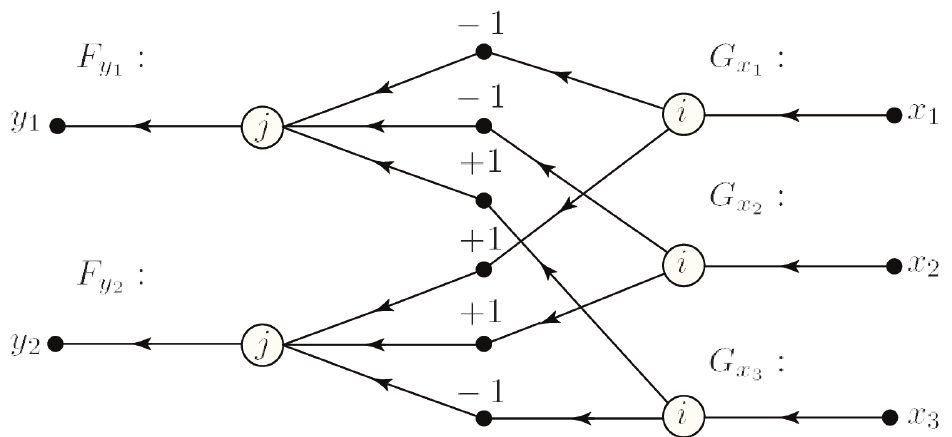}}
$$

Taking $X=Y=[1]$, we have the commutative, associative, unital monoid, with an involution, $({\mathbb Z}^{\otimes n})_{[1],[1]}$, its element are represented as $(F,G,\sigma,\varepsilon)$, $F,G$ rooted trees with 
$$
\mbox{labels $(F \backslash \partial F) \amalg (G \backslash \partial G) \to \{ 1,\ldots , n \}$, bijection $\sigma : \partial F \to \partial G$, and signs $\varepsilon : \partial F \to \{\pm 1\}$. }
$$
The relations are represented graphically as

\bigskip

\noindent (10.13) \quad $\delta$-{\bf unit}: 
$${\includegraphics[width=105mm]{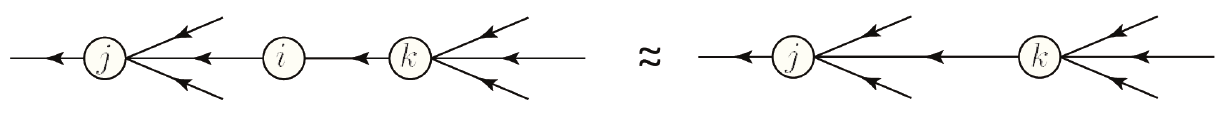}}$$

\bigskip

\noindent (10.14) \quad $\delta$-{\bf commutativity}: isomorphic data are equivalent.

\bigskip

\noindent (10.15) \quad $\delta$-{\bf associativity}: 
$$
{\includegraphics[width=70mm]{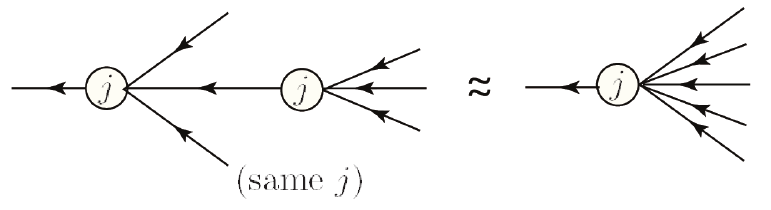}}
$$

\bigskip

\noindent (10.16) \quad {\bf cancellation}: 
$$
{\includegraphics[width=120mm]{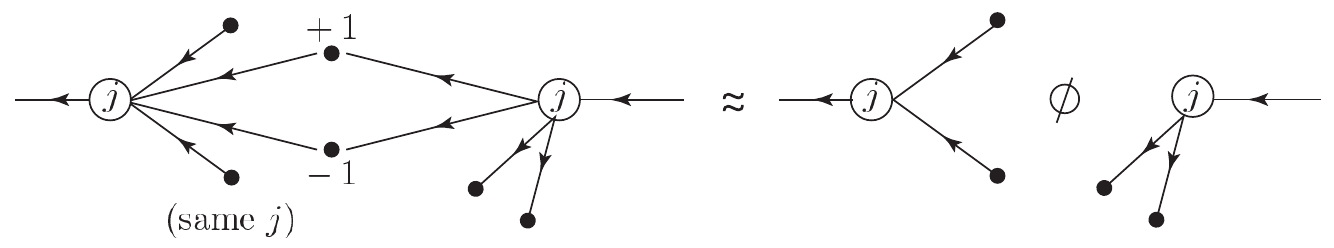}}
$$

\bigskip

\noindent (10.17) \quad {\bf commutativity}: 
$$
{\includegraphics[width=175mm]{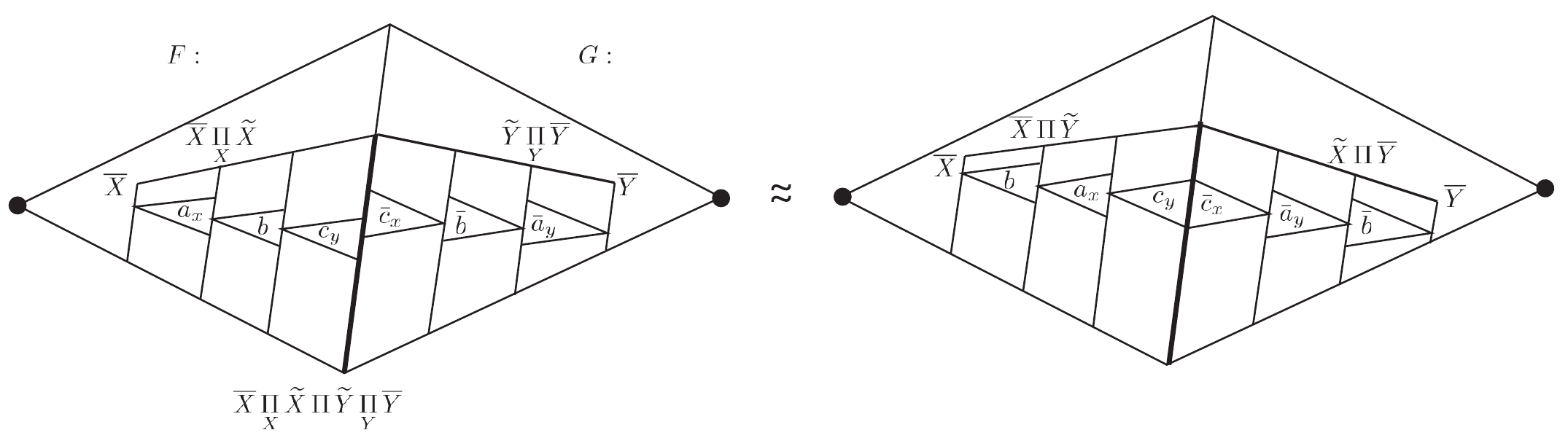}}
$$

Here $\overline X$ is a cut of the tree $F$, $\overline Y$ is a cut of the tree $G$, we have mappings
\setcounter{equation}{17}
\begin{equation}
\label{eq1018}
\xymatrix{
\underset{x \, \in \, X}{\coprod} \partial \, \bar c_x \equiv \overline X \ar[rd] \mbox{\hglue-1cm}&&\mbox{\hglue-1cm}\widetilde X \equiv \partial \, \bar b \equiv \underset{x \, \in \, X}{\coprod} \partial \, a_x \ar[ld] \\
&X
} , \xymatrix{
\underset{y \, \in \, Y}{\coprod} \partial \, \bar a_y \equiv \partial \, b \equiv \widetilde Y \ar[dr] \mbox{\hglue-1cm}&&\mbox{\hglue-1cm} \overline Y \equiv \underset{y \, \in \, Y}{\coprod} \partial \, c_y \ar[dl] \\
&Y
}
\end{equation}
and the bijection $\sigma$ identifies the part of $\partial F$ over $\overline X$ with the part of $\partial G$ over $\overline Y$, both being $\left( \overline X \, \underset{X}{\Pi} \, \widetilde X) \, \Pi \, (\widetilde Y \, \underset{Y}{\Pi} \, \overline Y \right)$.

\smallskip

The operation of multiplication is given by
\begin{equation}
\label{eq1019}
(F,G,\sigma,\varepsilon)_{/ \approx} \lhd \ (F',G',\sigma' , \varepsilon')_{/ \approx} \ = \left(F \underset{\partial F}{\lhd} F', G' \underset{\partial G'}{\lhd} G , \ \sigma \, \pi \, \sigma' , \ \varepsilon \cdot \varepsilon'\right)_{\!\!/ \approx}
\end{equation}
where $F \underset{\partial F}{\lhd} F'$ is the tree obtained by attaching a copy of $F'$ at every leaf $\in \partial F$ of $F$.

\smallskip

The involution is given by
\begin{equation}
\label{eq1020}
(F,G,\sigma,\varepsilon)_{/\approx}^t = (G,F,\sigma^{-1} , \varepsilon \circ \sigma^{-1})_{/\approx} \, .
\end{equation}
For each $i=1,2,\ldots ,n$, we can define an operation of addition on $({\mathbb Z}^{\otimes n})_{[1],[1]}$:
\begin{equation}
\label{eq1021}
(F_1 ,G_1 ,\sigma_1 , \varepsilon_1)_{/\approx} +_{(i)} (F_2 ,G_2 ,\sigma_2 , \varepsilon_2)_{/\approx} := (\delta_i \lhd (F_1,F_2) , \delta_i \lhd (G_1 , G_2) , \sigma_1 \amalg \sigma_2 , \varepsilon_1 \amalg \varepsilon_2)_{/\approx}  \, .
\end{equation}
With this operation of addition we obtain an ordinary commutative ring $({\mathbb Z}^{\otimes n})_{[1]}^{\{ i \}}$, with an involution. The symmetric group $S_n$ acts naturally on $({\mathbb Z}^{\otimes n})_{[1]}$, preserving multiplication and involution, and $g \in S_n$ takes $i$-th addition into $g(i)$-th addition:
\begin{equation}
\label{eq1022}
g \, (a_1 +_{(i)} a_2) = g(a_1) +_{(g(i))} g(a_2)
\end{equation}
so $g : ({\mathbb Z}^{\otimes n})_{[1]}^{\{ i \}} \to ({\mathbb Z}^{\otimes n})_{[1]}^{\{ g(i)\}}$ is a ring homomorphism.

\smallskip

Note that a prime ${\mathfrak p} \in {\rm spec} \, ({\mathbb Z}^{\otimes n})$ is a subset of $({\mathbb Z}^{\otimes n})_{[1]}$, which is a prime of all the rings $({\mathbb Z}^{\otimes n})_{[1]}^{\{i\}}$, $i =1,\ldots ,n$, and ${\mathfrak p}$ is closed with respect to the additions $+_{(i)}$ (as well as closed w.r.t. ``richer-avarages'').

\bigskip

In a similar way, we get from the rigs $(0,4)$, the ordinary commutative rings, with involution,
\begin{equation}
\label{eq1023}
\left( {\mathbb Z} \underset{\mathbb F}{\otimes} \{0,1\}\right)_{[1]} \subseteq \left( {\mathbb Z} \underset{\mathbb F}{\otimes} [0,1] \right)_{[1]} \subseteq \left( {\mathbb Z} \underset{\mathbb F}{\otimes} [0,\infty) \right)_{[1]} 
\end{equation}
and the continuous maps of compact Zariski topological spaces
\begin{equation}
\label{eq1024}
\xymatrix{
{\rm spec} \left( {\mathbb Z} \underset{\mathbb F}{\otimes} \{ 0,1 \} \right)_{[1]} &{\rm spec} \left( {\mathbb Z} \underset{\mathbb F}{\otimes} [ 0,1 ] \right)_{[1]} \ar[l] &{\rm spec} \left( {\mathbb Z} \underset{\mathbb F}{\otimes} [ 0,\infty ) \right)_{[1]} \ar[l] \\
\mbox{\begin{rotate}{90}$\subseteq$\end{rotate}} &\mbox{\begin{rotate}{90}$\subseteq$\end{rotate}} &\mbox{\begin{rotate}{90}$\subseteq$\end{rotate}} \\
{\rm spec} \, {\mathbb Z} \underset{\mathbb F}{\otimes} \{0,1\} &{\rm spec} \, {\mathbb Z} \underset{\mathbb F}{\otimes} [0,1] \ar[l] &{\rm spec} \, {\mathbb Z} \underset{\mathbb F}{\otimes} [0,\infty) \ar[l]
}
\end{equation}
(top line are spec's of ordinary $C$Ring, bottom line are spectra of $C{\mathcal G}{\mathbb R}$).

\smallskip

The group of positive reals ${\mathbb R}^+$ acts on the rig $[0,\infty)$ by automorphism
\begin{eqnarray}
\label{eq1025}
f : {\mathbb R}^+ &\to &{\rm Aut} \, ([0,\infty)) \\
\sigma &\mapsto &f^{\sigma} (x) = x^{\sigma} \nonumber
\end{eqnarray}
$$
f^{\sigma} (x_1 \cdot x_2) = f^{\sigma} (x_1) \cdot f^{\sigma} (x_2) \, , \quad f^{\sigma} (\max \, \{x_1 , x_2\}) = \max \, \{ f^{\sigma} (x_1) , f^{\sigma} (x_2)\}
$$
$$
f^{\sigma_1} (f^{\sigma_2} (x)) = f^{\sigma_1 \cdot \sigma_2} (x)
$$
and the ``fixed-field'' of this action is $\{0,1\}$. It follows that ${\mathbb R}^+$ acts on the compact spaces on the right of (\ref{eq1024}), preserving the fibers over the spaces of the left of (\ref{eq1024}).

\section{Completed vector bundles}

We fix $X = \{X_N , N \in I \, ; \ \pi_N^M : X_M \to X_N \, , \ M \geq N \} \in {\rm pro} \, C{\mathbb F} R^t$-Sch. For ${\mathcal U} \subseteq X_N$ open, let
\begin{equation}
\label{eq111}
S({\mathcal U}) = \left\{
\begin{matrix}
s=s^t \in {\mathcal O}_{X_N} ({\mathcal U})^+_{[1],[1]} \, , \ \forall \, M \geq N \, , \ \forall \, {\mathcal V} \subseteq (\pi_N^M)^{-1} \ {\mathcal U} \, , \ \forall \, a,a' \in {\mathcal O}_{X_M} ({\mathcal V})_{d,d'} : \\ { \ } \\
(\pi_N^M)^{\#} \, s \lhd a = (\pi_N^M)^{\#} \, s \lhd a' \Rightarrow a=a'
\end{matrix}
\right\} .
\end{equation}
These multiplicative sets give rise to a pre-sheaf ${\mathcal U} \mapsto S({\mathcal U})^{-1} \, {\mathcal O}_{X_N} ({\mathcal U})$ of $C{\mathbb F}R^t$ over $X_N$, and we let ${\mathcal K}_N$ denote the associated sheaf. We have an embedding ${\mathcal O}_{X_N} \hookrightarrow {\mathcal K}_N$, and we obtain an embedding of sheaves of groups over $X_N$
\begin{equation}
\label{eq112}
{\rm GL}_d ({\mathcal O}_{X_N}) \hookrightarrow {\rm GL}_d ({\mathcal K}_N) \, .
\end{equation}
The global sections of the quotient sheaf is the set of rank $d$ vector bundles on $X_N$ (trivialized at the generic points): 
\begin{equation}
\label{eq113}
D_d (X_N) = \Gamma (X_N , {\rm GL}_d ({\mathcal K}_N) / {\rm GL}_d ({\mathcal O}_{X_N})) \, .
\end{equation}
Its elements $D \in D_d (X_N)$ are represented as $D = \{ {\mathcal U}_{\alpha} , f_{\alpha} \}_{/\approx}$, where $X_N = \underset{\alpha}{\cup} \ {\mathcal U}_{\alpha}$ is an open cover, \qquad \qquad$f_{\alpha} \in {\rm GL}_d ({\mathcal K}_N)({\mathcal U}_{\alpha})$, and for all $\alpha , \beta : f_{\alpha}^{-1} \circ f_{\beta} \in {\rm GL}_d ({\mathcal O}_{X_N}) ({\mathcal U}_{\alpha} \cap {\mathcal U}_{\beta})$. We have the equivalence relation,

\bigskip

\noindent (11.4) \quad $\{{\mathcal U}_{\alpha} , f_{\alpha}\} \approx \{{\mathcal V}_{\beta} , g_{\beta} \}$ if and only if there is a common refinement $X_N = \underset{\gamma}{\cup} \ {\mathcal W}_{\gamma}$, and $u_{\gamma} \in {\rm GL}_d ({\mathcal O}_{X_N}) ({\mathcal W}_{\gamma})$, such that for ${\mathcal W}_{\gamma} \subseteq {\mathcal U}_{\alpha} \cap {\mathcal V}_{\beta}$, $f_{\alpha} = g_{\beta} \circ u_{\gamma} \in {\rm GL}_d ({\mathcal K}_N)({\mathcal W}_{\gamma})$.

\bigskip

We obtain the category of vector bundles over $X_N$, with objects $\underset{d}{\cup} \, D_d (X_N)$, and with arrows from $D = \{{\mathcal U}_{\alpha} , f_{\alpha} \}_{/\approx}$ to $D' = \{{\mathcal V}_{\beta} , g_{\beta} \}_{/\approx}$ given by
\begin{equation}
\setcounter{equation}{5}
\label{eq115}
{\rm Bun}_{X_N} (D,D') = \{ h \in {\mathcal K}_N (X_N)_{d',d} \, , \ \forall \, \alpha , \beta : g_{\beta}^{-1} \circ h \circ f_{\alpha} \in {\mathcal O}_{X_N} ({\mathcal V}_{\beta} \cap {\mathcal U}_{\alpha})_{d',d} \} \, .
\end{equation}
In particular, we get a partial order on $D_d (X_N)$,
\begin{equation}
\label{eq116}
D = \{ {\mathcal U}_{\alpha} , f_{\alpha} \}_{/\approx} \leq D' = \{ {\mathcal V}_{\beta} , g_{\beta}\}_{/\approx} \Leftrightarrow {\rm id}_d \in {\rm Bun}_{X_N} (D',D) \Leftrightarrow f_{\alpha}^{-1} \circ g_{\beta} \in {\mathcal O}_{X_N} ({\mathcal U}_{\alpha} \cap {\mathcal V}_{\beta})_{d,d'} \, .
\end{equation}
For $D = \{ {\mathcal U}_{\alpha} , f_{\alpha} \}_{/\approx} \in D_d (X_N)$, we have the subsheaf
\begin{equation}
\label{eq117}
{\mathcal O}_{X_N} (D)_{d'} = f_{\alpha} \circ ({\mathcal O}_{X_N} \vert_{{\mathcal U}_{\alpha}})_{d,d'} \subseteq ({\mathcal K}_N \vert_{{\mathcal U}_{\alpha}})_{d,d'}
\end{equation}
with ${\mathcal O}_{X_N} (D)_{d'} \circ ({\mathcal O}_{X_N})_{d',d''} \subseteq {\mathcal O}_{X_N} (D)_{d''}$.

\smallskip

Every $h \in {\rm Bun}_{X_N} (D,D')$ gives a mapping of sheaves
\begin{eqnarray}
\label{eq118}
{\mathcal O}_{X_N} (D)_{d''} &\to &{\mathcal O}_{X_N} (D')_{d''}  \\
s &\mapsto &h \circ s \nonumber
\end{eqnarray}
commuting with the right action of ${\mathcal O}_{X_N}$. We have
\begin{equation}
\label{eq119}
D \leq D' \Leftrightarrow {\mathcal O}_{X_N} (D)_{d'} \supseteq {\mathcal O}_{X_N} (D')_{d'} \, , \quad {\rm all} \ d' \, .
\end{equation}

For $M \geq N$ in $I$, we have pull-back of bundles
\begin{eqnarray}
\label{eq1110}
(\pi_N^M)^* : {\rm Bun}_{X_N} &\to &{\rm Bun}_{X_M} \\
(\pi_N^M)^* \{{\mathcal U}_{\alpha} , f_{\alpha}\}_{/\approx} &= &\{(\pi_N^M)^{-1} \, {\mathcal U}_{\alpha} , (\pi_N^M)^{\#} f_{\alpha} \}_{/\approx} \, . \nonumber
\end{eqnarray}
This is a functor, and an order preserving map $D_d (X_N) \to D_d (X_M)$. Since $I$ is directed, (9.1), we can easily pass to the co-limit category $\underset{N \in I}{\varinjlim} \, {\rm Bun}_{X_N}$, with rank-$d$ objects given by the partially ordered set $\underset{N \in I}{\varinjlim} \, D_d (X_N)$. But it seems that we should replace this co-limit by a kind of limit. For this we look at the set of monotone sequences
\begin{equation}
\label{eq1111}
{\mathcal B}_d (X) := \{D = \{ D_N\}_{N \geq N_0} \, , \ D_N \in D_d (X_N) \, , \ D_M \geq (\pi_N^M)^* \, D_N \quad {\rm for} \quad M \geq N \geq N_0 \}
\end{equation}
and its subset of bounded sequences
\begin{equation}
\label{eq1112}
{\mathcal B}_d^* (X) := \{D = \{ D_N\}_{N \geq N_0} \in {\mathcal B}_d (X) \, , \  \exists \, D_0 \in D_d (X_{N_0}) \, , \ D_M \leq (\pi_{N_0}^M)^* \, D_0 \quad {\rm for} \quad M \geq N_0 \} \, .
\end{equation}
For $N_1 \in I$, ${\mathcal U} \subseteq X_{N_1}$ open, and for $D = \{D_N\}_{N \geq N_0} \in {\mathcal B}_d(X)$, we have the ``$d_1$-sections of $D$ over ${\mathcal U}$'':
\begin{equation}
\label{eq1113}
{\mathcal O} (D,{\mathcal U})_{d_1} = \{ h \in {\mathcal K}_{N_1} ({\mathcal U})_{d,d_1} \, , \ (\pi_{N_1}^M)^* \, h \in {\mathcal O}_{X_M} (D_M)_{d_1} \quad {\rm for} \quad M \geq N_0 , N_1 \} \, .
\end{equation}
We define for $D = \{D_N \}_{N \geq N_0}$, $D' = \{D'_N\}_{N \geq N'_0} \in {\mathcal B}_d(X)$,
\begin{equation}
\label{eq1114}
D \leq D' \quad {\rm iff} \quad {\mathcal O} (D,{\mathcal U})_{d_1} \supseteq {\mathcal O} (D' , {\mathcal U})_{d_1} \quad \mbox{for all} \quad {\mathcal U} \subseteq X_{N_1} \ {\rm open}, \ {\rm all} \ d_1 \, .
\end{equation}
We get an equivalence relation on ${\mathcal B}_d (X)$
\begin{equation}
\label{eq1115}
D \approx D' \quad \mbox{iff} \quad D \leq D' \quad \mbox{and} \quad D' \leq D \, .
\end{equation}
We let ${\rm Bun}_d (X) = {\mathcal B}_d^* (X)_{/\approx}$ denote the equivalence classes, and we have an induce partial order $\leq$ on ${\rm Bund}_d (X)$. We have a category with objects $\underset{d}{\cup} \, {\rm Bun}_d (X)$, and arrows 
\begin{equation}
\label{eq1116}
{\rm Bun}_X (\{D_N\}_{/\approx} , \{D'_N\}_{/\approx}) = \left\{ h \in \varinjlim_{N \in I} {\mathcal K}_N(X_N)_{d',d} \, , \ h \circ {\mathcal O} (D,{\mathcal U})_{d_1} \subseteq {\mathcal O} (D' , {\mathcal U})_{d_1} \, , \quad {\rm all} \ {\mathcal U} \subseteq X_{N_1} \, {\rm open}, \ {\rm all} \ d_1 \right\} .
\end{equation}
Note that the group ${\rm GL}_d ({\mathcal K} (X)) = \underset{N \in I}{\varinjlim} \, {\rm GL}_d ({\mathcal K}_N (X_N))$ acts on ${\rm Bun}_d (X)$, and we can view ${\rm GL}_d ({\mathcal K} (X)) \backslash {\rm Bun}_d (X)$ as the isomorphism classes of rank-$d$ bundles.

\smallskip

For a number field $K$, and for $X = \{ X_N^K \} = \overline{{\rm spec} \, {\mathcal O}_K}$, we get
\begin{equation}
\label{eq1117}
{\rm Bun}_d \, (\overline{{\rm spec} \, {\mathcal O}_K}) \cong {\rm GL}_d ({\mathbb A}_K) / \prod_v {\rm GL}_d \, (\widehat{\mathcal O}_{K,v})
\end{equation}
where ${\mathbb A}_K$ is the ring of adeles of $K$, ${\rm GL}_d \, (\widehat{\mathcal O}_{K,v}) = \left\{ \begin{matrix} O(d) &v \ {\rm real} \hfill \\ U(d) &v \ {\rm complex} \end{matrix} \right.$, with the partial order
\begin{equation}
\label{eq1118}
g = (g_v) \leq g' = (g'_v) \Leftrightarrow g \cdot \widehat{\mathcal O}_{K,v}^{\oplus d} \supseteq g' \cdot \widehat{\mathcal O}_{K,v}^{\oplus d} \quad \mbox{for all $v$},
\end{equation}
with $\widehat{\mathcal O}_{K,v}^{\oplus d} = \left\{ (x_1 , \ldots , x_d \} \in \widehat K_v^{\oplus d} \, , \ \underset{i=1}{\overset{d}\sum} \, \vert x_i \vert^2 \leq 1 \right\}$ for $v$ real or complex, and with its natural 

\smallskip

\noindent ${\rm GL}_d (K) = \underset{N}{\varinjlim} \, {\rm GL}_d ({\mathcal K}_N (X_N))$ action (note that ${\mathcal K}_N (X_N) \equiv K$, all $N$). 

\smallskip

E.g. For $d=1$ we get the Picard group of isomorphism classes of line bundles
\begin{equation}
\label{eq1119}
{\rm Pic} \left( \overline{{\rm spec} \, {\mathcal O}_K} \right) = K^* \backslash {\mathbb A}_K^* / \prod_v \widehat{\mathcal O}_{K,v}^* \overset{\rm deg}{-\!\!\!-\!\!\!-\!\!\!-\!\!\!-\!\!\!\twoheadrightarrow} {\rm Pic} \, \left(\overline{{\rm spec} \, {\mathbb Z}}\right) = {\mathbb R}^+ \, .
\end{equation}
The kernel ${\rm Pic}^{\circ} \left( \overline{{\rm spec} \, {\mathcal O}_K} \right) = \ker \deg$, is a compact group, an extension of the ideal class group by the Dirichlet units torus:
\begin{equation}
\label{eq1120}
* \to \mu_K \to {\mathcal O}_K^* \to \left[ \prod_{ {\overset{v \, {\rm real}}{{\rm or \,
complex}}} } {\mathbb R}^+ \right]^{\circ} \to {\rm Pic}^{\circ} \left( \overline{{\rm spec} \, {\mathcal O}_K} \right) \to {\rm class}_K \to * \ .
\end{equation}

\end{document}